\documentclass[12pt,a4paper]{amsart}

\usepackage{soul}
\usepackage{lscape}

\usepackage{DALdefams}

\newcommand {\fBj}{{\mathfrak{Bj}}}
\newcommand {\fbj}{{\mathfrak{bj}}}
\newcommand {\fbr}{{\mathfrak{br}}}
\newcommand {\fbrj}{{\mathfrak{brj}}}
\newcommand {\fBRJ}{{\mathfrak{BRJ}}}
\newcommand {\fel}{{\mathfrak{el}}}
\newcommand {\fMe}{{\mathfrak{Me}}}

\begin{document}

\title[New simple Lie superalgebras]
{New simple modular Lie superalgebras as generalized prolongs}

\author{Sofiane Bouarroudj${}^1$, Pavel Grozman${}^2$, Dimitry Leites${}^3$}

\address{${}^1$Department of Mathematics, United Arab Emirates University, Al
Ain, PO. Box: 17551; Bouarroudj.sofiane@uaeu.ac.ae\\ ${}^2$Equa
Simulation AB, Stockholm, Sweden; pavel@rixtele.com\\
${}^3$MPIMiS, Inselstr. 22, DE-04103 Leipzig, Germany\\ on leave
from Department of Mathematics, University of Stockholm, Roslagsv.
101, Kr\"aft\-riket hus 6, SE-106 91 Stockholm, Sweden;
mleites@math.su.se, leites@mis.mpg.de}



\begin{abstract} Over algebraically closed fields of
characteristic $p>2$, prolongations of the simple finite dimensional
Lie algebras and Lie superalgebras with Cartan matrix are studied
for certain simplest gradings of these algebras. Several new simple
Lie superalgebras are discovered, serial and exceptional, including
superBrown and superMelikyan superalgebras. Simple Lie superalgebras
with Cartan matrix of rank 2 are classified.
\end{abstract}

\thanks{We are thankful to I.~Shchepochkina for help;
DL is thankful to MPIMiS, Leipzig, for financial support and most
creative environment.}


\maketitle

\markboth{\itshape Sofiane Bouarroudj\textup{,} Pavel
Grozman\textup{,} Dimitry Leites}{{\itshape New simple modular Lie
superalgebras}}

\thispagestyle{empty}

\section{Introduction}

\ssec{Setting} We use standard notations of \cite{FH, S}; for the
precise definition (algorithm) of generalized
Cartan-Tanaka--Shchepochkina (CTS) complete and partial
prolongations, and algorithms of their construction, see
\cite{Shch}. Hereafter $\Kee$ is an algebraically closed field of
characteristic $p>2$, unless specified. Let $\fg'=[\fg, \fg]$, and
$\fc(\fg)=\fg\oplus\mathfrak{center}$, where
$\dim\mathfrak{center}=1$. Let ${}^{n)}\fg$ denote the incarnation
of the Lie (super)algebra $\fg$ with the $n)$th Cartan matrix, cf.
\cite{GL4, BGL1}. On classification of simple vectorial Lie
superalgebras with polynomial coefficients (in what follows referred
to as {\it vectorial Lie superalgebras of polynomial vector fields}
over $\Cee$, see \cite{LSh, K3}).

The works of S.~Lie, Killing and \`E.~Cartan, now classical,
completed classification  over $\Cee$ of
\begin{equation}
\label{1} \text{{\bf simple Lie algebras of finite dimension and
of polynomial vector fields}.}
\end{equation}

Lie algebras and Lie superalgebras over fields in characteristic
$p>0$, a.k.a. {\it modular} Lie (super)algebras, were first
recognized and defined in topology, in the 1930s. The {\bf simple}
Lie algebras drew attention (over finite fields $\Kee$) as a step
towards classification of simple finite groups, cf. \cite{St}. Lie
{\it super}algebras, even simple ones  and even over $\Cee$ or
$\Ree$, did not draw much attention of mathematicians until their
(outstanding) usefulness was observed by physicists in the 1970s.
Meanwhile mathematicians kept discovering new and new examples of
simple modular Lie algebras until Kostrikin and Shafarevich
(\cite{KS}) formulated a conjecture embracing all previously found
examples for $p>7$. Its generalization reads: {\sl select a
$\Zee$-form $\fg_\Zee$ of every $\fg$ of type\footnote{Observe that
the algebra of divided powers (the analog of the polynomial algebra
for $p>0$) {\it and hence all prolongs} (Lie algebras of vector
fields) acquire one more --- shearing --- parameter:
$\underline{N}$, see \cite{S}.} $(\ref{1})$, take
$\fg_\Kee:=\fg_\Zee\otimes_\Zee\Kee$ and its simple finite
dimensional subquotient $\fs\fii(\fg_\Kee)$ (there can be several
such $\fs\fii(\fg_\Kee)$). Together with deformations\footnote{It is
not clear, actually, if the conventional notion of deformation can
always be applied if $p>0$ (for the arguments, see \cite{LL}; cf.
\cite{Vi}); to give the correct (better say, universal) notion is an
open problem, but in some cases it is applicable, see \cite{BGL4}.}
of these examples we get in this way all simple finite dimensional
Lie algebras over algebraically closed fields if $p>5$. If $p=5$, we
should add to the above list Melikyan's examples}.

Having built upon ca 30 years of work of several teams of
researchers, and having added new ideas and lots of effort, Block,
Wilson, Premet and Strade proved the generalized KSh conjecture for
$p>3$, see \cite{S}. For $p\leq 5$, the above KSh-procedure does not
produce all simple finite dimensional Lie algebras; there are other
examples. In \cite{GL4}, we returned to \'E.~Cartan's description of
$\Zee$-graded Lie algebras as CTS prolongs, i.e., as subalgebras of
vectorial Lie algebras preserving certain distributions; we thus
interpreted the \lq\lq mysterious" at that moment exceptional
examples of simple Lie algebras for $p=3$ (the Brown, Frank,
Ermolaev and Skryabin algebras), further elucidated Kuznetsov's
interpretation \cite{Ku1} of Melikyan's algebras (as prolongs of the
nonpositive part of the Lie algebra $\fg(2)$ in one of its
$\Zee$-gradings) and discovered three new series of simple Lie
algebras. In \cite{BjL}, the same approach yielded $\fbj$, a simple
super versions of $\fg(2)$, and $\fBj(1; N|7)$, a simple $p=3$ super
Melikyan algebra. Both $\fbj$ and $\fBj(1; N|7)$ are indigenous to
$p=3$, the case where $\fg(2)$ is not simple.

\ssec{Classification: Conjectures and results} Recently, Elduque
considered super analogs of the exceptional simple Lie algebras; his
method leads to a discovery of 10 new simple (presumably,
exceptional) Lie superalgebras for $p=3$. For a description of the
Elduque superalgebras, see \cite{CE, El1, CE2, El2}; for their
description in terms of Cartan matrices and analogs of Chevalley
relations and notations we use in what follows, see \cite{BGL1,
BGL2}.

In \cite{L}, a super analog of the KSh conjecture embracing all
types of simple (finite dimensional) Lie superalgebras is formulated
based on an entirely different idea in which the CTS prolongs play
the main role:

\so{For every simple finite dimensional Lie (super)algebra of the
form $\fg(A)$, take its non-positive part with respect to a certain
simplest $\Zee$-grading, consider its complete and partial prolongs
and take their simple subquotients.}

The new examples of simple modular Lie superalgebras ($\fBRJ$,
$\fBj(3;\underline{N}|3)$, $\fBj(3;\underline{N}|5)$) support this
conjecture. (This is how Cartan got all simple $\Zee$-graded Lie
algebras of polynomial growth and finite depth --- the Lie algebras
of type $(\ref{1})$ --- at the time when the root technique was not
discovered yet.)

\subsubsection{Yamaguchi's theorem (\cite{Y})}  This theorem,
reproduced in \cite{GL4, BjL}, states that {\sl for almost all
simple finite dimensional Lie algebras $\fg$ over $\Cee$ and their
$\Zee$-gradings $\fg=\mathop{\oplus}\limits_{-d\leq i}\fg_i$ of
finite depth $d$, the CTS prolong of
$\fg_{\leq}=\mathop{\oplus}\limits_{-d\leq i\leq 0}\fg_i$ is
isomorphic to $\fg$, the rare exceptions being two of the four
series of simple vectorial algebras} (the other two series being
partial prolongs).

\subsubsection{Conjecture} In the following theorems, we
present the results of \texttt{SuperLie}-assisted (\cite{Gr})
computations of the CTS-prolongs of the non-positive parts of the
simple finite dimensional Lie algebras and Lie superalgebras
$\fg(A)$; we have only considered $\Zee$-grading corresponding to
each (or, for larger ranks, even certain {\it selected}) of the {\bf
simplest} gradings $r=(r_1, \dots, r_{\rk\fg})$, where all but one
coordinates of $r$ are equal to 0 and only one --- {\it selected}
--- is equal to 1, and where we set $\deg X_i^{\pm}=\pm r_i$ for
the Chevalley generators $X_i^{\pm}$ of $\fg(A)$, see \cite{BGL1}.

\so{Other gradings (as well as algebras $\fg(A)$ of higher ranks) do
not yield new simple Lie (super)algebras as prolongs of the
non-positive parts}.

\ssbegin{Theorem}\label{th3} The CTS prolong of the nonpositive part
of $\fg$ returns $\fg$ in the following cases: $p=3$ and
$\fg=\ff(4)$, $\fe(6)$, $\fe(7)$ and $\fe(8)$ considered with the
$\Zee$-grading with one selected root corresponding to the endpoint
of the Dynkin diagram.
\end{Theorem}

\subsubsection{Conjecture}[The computer got stuck here, after weeks of
computations]\label{cj1} To the cases of Theorem \ref{th3}, one can
add the case for $p=5$ and $\fg=\fel(5)$ {\em (see \cite{BGL2})} in
its $\Zee$-grading with only one odd simple root and with one
selected root corresponding to any endpoint of the Dynkin diagram.

\ssbegin{Theorem}\label{th4} Let $p=3$. For the previously known (we
found more, see Theorems \ref{th6}, \ref{th9}) simple finite
dimensional Lie superalgebras $\fg$ of rank $\leq 3$ with Cartan
matrix and for their simplest gradings $r$, the CTS prolongs (of the
non-positive part of $\fg$) different from $\fg$ are given in the
following table elucidated below.
\end{Theorem}

\ssec{Melikyan superalgebras for $p=3$} There are known the two
constructions of the Melikyan algebra $\fMe(5;
\underline{N})=\mathop{\oplus}\limits_{i\geq -2}\fMe(5;
\underline{N})_i$, defined for $p=5$:

1) as the CTS  prolong of the triple $\fMe_0=\fc\fvect(1;
\underline{1})$, $\fMe_{-1}=\cO(1;\underline{1})/\text{const}$ and
the trivial module $\fMe_{-2}$, see \cite{S}; this construction
would be a counterexample to our conjecture were there no
alternative:

2) as the {\it complete} CTS prolong of the non-negative part of
$\fg(2)$ in its grading $r=(01)$, with $\fg(2)$ obtained now as a
{\it partial} prolong, see \cite{Ku1, GL4}.

In \cite{BjL}, we have singled out $\fBj(1; \underline{N}|7)$ as a
$p=3$ simple analog of $\fMe(5; \underline{N})$ as a partial CTS
prolongs of the pair (the negative part of $\fk(1;
\underline{N}|7)$, $\fBj(1; \underline{N}|7)_0=\fp\fgl(3)$), and
$\fbj$ as a $p=3$ simple analog of $\fg(2)$ whose non-positive part
is the same as that of $\fBj(1; \underline{N}|7)$, i.e., $\fbj$ and
$\fBj(1; \underline{N}|7)$ are analogs of the construction 2).

The original Melikyan's construction 1) also has its super analog
for $p=3$ (only in the situation described in Theorem \ref{th6}) and
it yields a new series of simple Lie superalgebras as the complete
prolongs, with another simple analog of $\fg(2)$ as a partial
prolong.

Recall (\cite{BGL1}) that we normalize the Cartan matrix $A$ so that
$A_{ii}=1$ or $0$  if the $i$th root is odd,  whereas if the $i$th
root is even, we set $A_{ii}=2$ or $0$ in which case we write $\ev$
instead of $0$ in order not to confuse with the case of odd roots.

\ssbegin{Theorem}\label{th6} A $p=3$ analog of the construction $1)$
of the Melikyan algebra is given by setting $\fg_0=\fc\fk(1;
\underline{1}|1)$, $\fg_{-1}=\cO(1;\underline{1}|1)/\text{const}$
and $\fg_{-2}$ being the trivial module. It yields a simple super
Melikyan algebra that we denote by $\fMe(3; \underline{N}|3)$,
non-isomorphic to a superMelikyan algebra $\fBj(1;
\underline{N}|7)$.

The partial prolong of the non-positive part of $\fMe(3;
\underline{N}|3)$ is a new (exceptional) simple Lie superalgebra
that we denote by $\fbrj(2;3)$. This $\fbrj(2;3)$ has the three
Cartan matrices:
$1)\ \begin{pmatrix}0&-1\\
-2&1\end{pmatrix}$ and $2)\ \begin{pmatrix}0&-1\\
-1&\ev\end{pmatrix}$ joined by an odd reflection, and $\begin{pmatrix}1&-1\\
-1&\ev\end{pmatrix}$.  It is a super analog of the Brown algebra
$\fbr(2)=\fbrj(2;3)_\ev$, its even part.

The CTS prolongs for the simplest gradings $r$ of
${}^{1)}\fbrj(2;3)$ returns known simple Lie superalgebras, whereas
the CTS prolong for a simplest grading $r$ of ${}^{2)}\fbrj(2;3)$
returns, as a partial prolong, a new simple Lie superalgebra we
denote $\fBRJ$.

Unlike $\fbr(2)$, the Lie superalgebra $\fbrj(2;3)$ has analogs for
$p\neq 3$, e.g., for $p=5$, we get a new simple Lie superalgebra
$\fbrj(2;5)$ such that $\fbrj(2;5)_\ev=\fsp(4)$ with the two Cartan
matrices
$1)\ \begin{pmatrix}0&-1\\
-2&1\end{pmatrix}$ and $2)\ \begin{pmatrix}0&-4\\
-3&2\end{pmatrix}$. The CTS prolongs of $\fbrj(2;5)$ for all its
Cartan matrices and the simplest $r$ return $\fbrj(2;5)$.
\end{Theorem}

Having got this far, it was impossible not to try to get
classification of simple $\fg(A)$'s. Here is its beginning part, see
\cite{BGL5}.

\ssbegin{Theorem}\label{th9} If $p>5$, every finite dimensional
simple Lie {\em super}algebra with a $2\times 2$ Cartan matrix is
isomorphic to $\fosp(1|4)$, $\fosp(3|2)$, or $\fsl(1|2)$. If $p=5$,
we should add $\fbrj(2;5)$. If $p=3$, we should add $\fbrj(2;3)$.
\end{Theorem}

\begin{Remark} For details of description  of the new simple Lie superalgebras
of types $\fBj$ and $\fMe$ and their subalgebras, in particular,
presentations of $\fbrj(2;3)$ and $\fbrj(2;5)$, and proof of Theorem
\ref{th9} and its generalization for higher ranks,  see \cite{BGL4,
BGL5}.
\end{Remark}

The new simple Lie superalgebras obtained are described in the next
subsections. \footnotesize
$$
\renewcommand{\arraystretch}{1.4}
\begin{tabular}{|l|l|l|l|}
\hline $\fg$&Cartan matrix&$r$&prolong\\ \hline
$\fosp(3|2)$&$\begin{matrix}\begin{pmatrix} 0 & -1  \\ -2 & 2
\end{pmatrix}\\
\begin{pmatrix} 0 & -1  \\ -1 &1
\end{pmatrix}\end{matrix}$&$\begin{matrix}(10)\\(01)\\(10)\\(01)\end{matrix}$&
$\begin{matrix}\fk(1|3)\\ \fk(1|3;1)\\ \fosp(3|2)\\
\fk(1|3;1)\end{matrix}$\\ \hline $\fsl(1|2)$&$\begin{pmatrix} 0 &
-1  \\ -1 & 2
\end{pmatrix}$&$\begin{matrix}(10)\\(01)\end{matrix}$&$\begin{matrix}\fvect(0|2)\\
\fvect(1|1)\end{matrix}$\\ &$\begin{pmatrix} 0 & 1  \\ -1 & 0
\end{pmatrix}$&$\begin{matrix}(10)\\(01)\end{matrix}$&$\fvect(1|1)$\\
\hline $\fosp(1|4)$&$\begin{pmatrix} 2 & -1  \\ -1 & 1
\end{pmatrix}$&$\begin{matrix}(10)\\(01)\end{matrix}$&$\begin{matrix}\fk(3|1)\\
\fosp(1|4)\end{matrix}$\\ \hline
$\fbrj(2;3)$&$\begin{matrix}\begin{pmatrix} 0 & -1  \\ -2 & 1
\end{pmatrix}\\
\begin{pmatrix} \ev & -1  \\ -1 &0
\end{pmatrix}\end{matrix}$&$\begin{matrix}(10)\\(01)\\(10)\\(01)\end{matrix}$&
$\begin{matrix}\mathfrak{Me}(3;N|3)\\ \mathfrak{Brj}(4|3)\\ \mathfrak{Brj}(4; \underline{N}|3) \\
\mathfrak{Brj}(3; \underline{N}|4)\supset\fBRJ
 \end{matrix}$\\
\hline $\fbrj(2;3)$&$\begin{matrix}\begin{pmatrix} \ev & -1\\
-1 & 1
\end{pmatrix}
\end{matrix}$&$\begin{matrix}(10)\\(01)\end{matrix}$&
$\begin{matrix}\mathfrak{Brj}(3; \underline{N}|3)\\
\mathfrak{Brj}(3; \underline{N}|4)\supset\fBRJ\end{matrix}$\\
\hline $\fbrj(2;5)$&$\begin{matrix}\begin{pmatrix} 0 & -1  \\ -2 & 1
\end{pmatrix}\\
\begin{pmatrix} 0 & 1  \\ -3 &2
\end{pmatrix}\end{matrix}$&$\begin{matrix}(10)\\(01)\\(10)\\(01)\end{matrix}$&
$\begin{matrix}\mathfrak{brj}(2;5)\\ \mathfrak{brj}(2;5)\\ \mathfrak{brj}(2;5)\\
\mathfrak{brj}(2;5)\end{matrix}$\\ \hline
\end{tabular}
$$

$$
\renewcommand{\arraystretch}{1.4}
\begin{tabular}{|l|l|l|l|}
\hline $\fg$&Cartan matrix&$r$&prolong\\ \hline
$\fsl(1|3)$&$\begin{pmatrix} 0 & -1 & 0\\ -1 & 2&-1\\0&-1&2
\end{pmatrix}$&$\begin{matrix}(100)\\(010)\\(001)\end{matrix}$&$\begin{matrix}
\fvect(0|3)\\\fsl(1|3)\\ \fvect(2|1)\end{matrix}$\\
&$\begin{pmatrix} 0 & -1  &0\\ -1 & 0&-2\\0&-1&2
\end{pmatrix}$&$\begin{matrix}(100)\\(010)\\(001)\end{matrix}$&$\begin{matrix}\fvect(2|1)\\
\fsl(1|3)\\ \fvect(2|1)\end{matrix}$\\ \hline $\fpsl(2|2)$&any
matrix&$\begin{matrix}(100)\\(010)\\(001)\end{matrix}$&
$\begin{matrix}\fsvect(1|2)\\ \fpsl(2|2)\\
\fsvect(1|2)\end{matrix}$\\
\hline $\fosp(1|6)$&$\begin{pmatrix} 2 & -1 & 0\\ -1 & 2&-1\\0&-1&1
\end{pmatrix}$&$\begin{matrix}(100)\\(010)\\(001)\end{matrix}$&$\begin{matrix}\fk(5|1)\\
\fosp(1|6)\\ \fosp(1|6)\end{matrix}$\\ \hline $\fosp(3|4)$&$
\left.\begin{matrix} \begin{pmatrix}2 & -1 & 0\\ -1 & 0&-1\\0&-2&2
\end{pmatrix}\\
\begin{pmatrix}
0&-1&0\\ -1&0&1\\ 0&-1&1
\end{pmatrix}\end{matrix}\right \}$&$\begin{matrix}(100)\\
(010)\\ (001)\end{matrix}$&$\begin{matrix}\fk(3|3)\\ \fosp(3|4)\\
\fosp(3|4)\end{matrix}$\\ &$\begin{pmatrix} 0 & -1 & 0\\ -1 &
2&-1\\0&-1&1
\end{pmatrix}$&$\begin{matrix}(100)\\
(010)\\ (001)\end{matrix}$&$\fosp(3|4)$\\
\hline $\fosp(5|2)$&$ \left.\begin{matrix} \begin{pmatrix} 2 & -1 &
0\\ -1 & 0&1\\0&-1&1
\end{pmatrix}\\
\begin{pmatrix}
0&-1&0\\ -1&0&1\\ 0&-2&2
\end{pmatrix}\end{matrix}\right \}$&$\begin{matrix}(100)\\
(010)\\ (001)\end{matrix}$&$\fosp(5|2)$\\ &$\begin{pmatrix} 0 & -1
& 0\\ -1 & 2&-1\\0&-2&2
\end{pmatrix}$&$\begin{matrix}\left .\begin{matrix}(100)\\
(010)\end{matrix}\right \}\\
(001)\end{matrix}$&$\begin{matrix}\fosp(5|2)\\
\fk(1|5)\end{matrix}$\\
\hline $ \begin{matrix}\fosp(4|2;\alpha)\\
\alpha\;\text{generic}\end{matrix}$&$
\begin{matrix}1)\;\begin{pmatrix}
 2 & -1 & 0 \\ \alpha& 0 & -1-\alpha \\ 0 & -1 & 2
\end{pmatrix}\\
2)\;\begin{pmatrix} 0 & 1 & -1-\alpha \\ -1 & 0 & -\alpha \\
-1-\alpha & \alpha & 0
\end{pmatrix}\end{matrix}$&$\begin{matrix}(100)\\
(010)\\ (001)\end{matrix}$&$\fosp(4|2;\alpha)$\\
\hline
$\begin{matrix} \fosp(4|2;\alpha)\\
\alpha=0,-1
\end{matrix}$&$ \begin{array}{l} 1) \text{ The simple part of
${}^{1)}\fosp(4|2;\alpha)$ is } \fsl(2|2);\\
\text{for the CTS of $\fpsl(2|2)$, see
above }\\
2)\; {}^{2)}\fosp(4|2;\alpha)\simeq\fsl(2|2); \\
\text{for the CTS of $\fsl(2|2)$, see
above }\\
\end{array}$&&\\
\hline
\end{tabular}
$$

$$
\renewcommand{\arraystretch}{1.4}
\begin{tabular}{|l|l|l|l|}
\hline $\fg$&Cartan matrix&$r$&prolong\\
\hline
$\fosp(2|4)$&$
\begin{matrix}1)\;\begin{pmatrix} 0 & 1 & 0 \\ -1 & 2 & -2 \\ 0 & -1 &2
\end{pmatrix}\\
2)\; \begin{pmatrix} 0 & 1 & 0 \\ -1 & 0 & 2 \\ 0 & -1 & 2
\end{pmatrix}\\
3)\; \begin{pmatrix} 0 & -2 & 1 \\ -2 & 0 & 1 \\ -1 & -1 & 2
\end{pmatrix}\end{matrix}$&$\begin{matrix}\left. \begin{matrix}(100)\\
(010)\\ (001)\end{matrix}\right\}\\ \left. \begin{matrix}(100)\\
(010)\\ (001)\end{matrix}\right\}\\ \left. \begin{matrix}(100)\\
(010)\\ (001)\end{matrix}\right\}\end{matrix}$&$
\begin{matrix}\begin{cases}\fosp(2|4)&\\
\begin{cases}\fosp(2|4)&\text{if $p>3$}\\
\fBj(3; N|3)&\text{if $p=3$}
\end{cases}&\\
\fosp(2|4)&\\
\end{cases}\\
\begin{cases}\fk(3|2)&\\
\begin{cases}\fosp(2|4)&\text{if $p>3$}\\
\fBj(3; N|3)&\text{if $p=3$}
\end{cases}\\
\fosp(2|4)&\\
\end{cases}\\
\begin{cases}\fosp(2|4)&\\
\fosp(2|4)&\\ \fk(3|2)&\\
\end{cases}
\end{matrix}$\\
\hline $\fg(2|3)$&$3)\; \begin{pmatrix}
    0&0&-1 \\
    0&0&-2 \\
    -1&-2&2
  \end{pmatrix}$&$\begin{matrix}(100)\\
  (010)\\
(001)\end{matrix}$& $\begin{matrix}\fBj(2|4)\\
\fBj(3|5)\\ \fbj\end{matrix}$\\ \hline
\end{tabular}
$$

\normalsize

\ssec{A description of $\fBj(3; N|3)$} For $\fg={}^{1)}\fosp(2|4)$
and $r=(0,1,0)$, we have the following realization of the
non-positive part:
\begin{equation}\label{nonpos}
\small
\renewcommand{\arraystretch}{1.4}
\begin{tabular}{|l|l|} \hline
$\fg_{i}$&the generators (even $\mid$ odd) \\ \hline \hline
$\fg_{-2}$&$Y_{6}=\partial_1\;\mid\; Y_{8}=\partial_4$\\ \hline
$\fg_{-1}$&$Y_2=\partial_2,\; Y_{5}=x_{2}\partial_{1} +
x_{5}\partial_{4} + \partial_{3},\;  \;\mid\; Y_4=\partial_5,\;
Y_{7}=2\, x_{2}\partial_{4} + \partial_{6},$\\ \hline $\fg_{0}\simeq
$& $Y_3={{x_{2}}^2}\partial_{1}
      +  x_{2} x_{5} \partial_{4}
      +  x_{2}\partial_{3} + 2x_{5}\partial_{6},\;
      Z_3={{x_{3}}^2}\partial_{1}
      +  2x_{3} x_{6}\partial_{4}
      +  x_{3}\partial_{2} + 2x_{6}\partial_{5}$\\
$\fsl(1|1)\oplus \fsl(2)\oplus \Kee$&$H_2=2\, x_{1}\partial_{1} +
2\, x_{2}\partial_{2}
   + x_{4}\partial_{4} + x_{5}\partial_{5}
   + 2\, x_{6}\partial_{6},\; H_1=[Z_1, Y_1], H_3=[Z_3,Y_3] \;\mid\;$\\
   & $ Y_1= x_{1}\partial_{4} + 2\, x_{2}\partial_{5}
   + x_{3}\partial_{6},\; Z_1=
  2\, x_{4}\partial_{1} + 2\, x_{5}\partial_{2}
   + x_{6}\partial_{3} $
   \\
\hline
\end{tabular}
\end{equation}
The $\fg_0$-module $\fg_{-1}$ is irreducible, having one highest
weight vector $Y_2$.

Let $p=3$. The CTS prolong gives $\sdim(\fg_1)=4|4$. The
$\fg_0$-module $\fg_1$ has the following two lowest weight
vectors: $$
\small
\renewcommand{\arraystretch}{1.4}
\begin{tabular}{|l|l|} \hline
$V_1'$ & $ x_{1} x_{2}  \partial_{4} + 2\, x_{1}\partial_{6} + 2\,
{{x_{2}}^2}  \partial_{5}+ x_{2}  x_{3}  \partial_{6}  $\\
   $V_1''$ & $  x_{1}x_{2}  \partial_{1}
           +  x_{1}  x_{5} \partial_{4}+ 2\, x_{2}  x_{4} \partial_{4}
   + x_{1}\partial_{3}+ {{x_{2}}^2} \partial_{2}
   + 2\, x_{2}  x_{5} \partial_{5} + x_{2}  x_{6} \partial_{6}
   +x_{3}  x_{5} \partial_{6}+ x_{4}\partial_{6}$\\
\hline
\end{tabular}
$$ Since $\fg_1$ generates the positive part of the CTS prolong,
$[\fg_{-1},\fg_1]=\fg_0,$ and $[\fg_{-1},\fg_{-1}]=\fg_{-2}$, the
standard criteria of simplicity ensures that the CTS prolong is
simple. Since none of the $\Zee$-graded Lie superalgebras over
$\Cee$ of polynomial growth and finite depth has grading of this
form (with $\fg_0\simeq\fsl(1|1)\oplus \fsl(2)\oplus \Kee$), we
conclude that this Lie superalgebra is new. We denote it by
$\fBj(3;N|3)$, where $N$ is the shearing parameter of the even
indeterminates. Our calculations show that $N_2=N_3=1$ always. For
$N_1=1$, 2, the super dimensions of the positive components of
$\fBj(3;N|3)$ are given in the following tables:
$$
\renewcommand{\arraystretch}{1.4}
\begin{tabular}{|l|c|c|c|c|c|c|c|c|c|} \hline
$N=(1,1,1)$ &$\fg_1$ & $\fg_2$& $\fg_3$&$\fg_4$ & $\fg_5$ &
$\fg_6$&--&--&--\\ \hline $\sdim$ &$4|4$ & $5|5$& $4|4$& $4|4$&
$2|2$& $0|3$&&&\\ \hline
 \hline
$N=(2,1,1)$ &$\fg_1$ & $\fg_2$& $\fg_3$&$\fg_4$ & $\fg_5$ &
$\fg_6$&$\cdots$ & $\fg_{11}$ & $\fg_{12}$\\ \hline $\sdim$ &$4|4$
& $5|5$& $4|4$& $5|5$& $4|4$ & $5|5$ &$\cdots$& $2|2$ & $0|3$\\
\hline
\end{tabular}
$$ Let $V_i'$, $V_i''$ and $V_i'''$ be the lowest height vectors
of $\fg_i$ with respect to $\fg_0$. For $N=(1,1,1)$, these vectors
are as follows: $$ \footnotesize
\renewcommand{\arraystretch}{1.4}
\begin{tabular}{|l|l|} \hline
$\fg_i$ & lowest weight vectors\\ \hline $V_2'$ & $  {{x_{1}}^2}
\partial_{4}
        + 2\, x_{1}  x_{2}  \partial_{5} + x_{1}  x_{3}  \partial_{6}
        +  x_{2} {{x_{3}}^2}  \partial_{6}$\\
$V_2''$ &$ x_{1} {{x_{2}}^2}  \partial_{1}
 + x_{1} x_{2}x_{5}  \partial_{4}
  +  x_{1}x_{2} \partial_{3}
+ 2\, x_{1}  x_{5} \partial_{6}
 + {{x_{2}}^2} x_{3}  \partial_{3}
  + 2\, {{x_{2}}^2}  x_{5}  \partial_{5}
   +  x_{2} x_{3}  x_{5}  \partial_{6} $\\
$V_2'''$ & $ {{x_{1}}^2} \partial_{1}
   +{{x_{2}}^2}  {{x_{3}}^2}  \partial_{1}
   +  {{x_{2}}^2} x_{3}  \partial_{2}
   + 2\, {{x_{2}}^2}  x_{6}  \partial_{5}
   + 2\, x_{1}  x_{2}  \partial_{2} + 2\,
 x_{1}  x_{3}\partial_{3} + 2\,  x_{1}  x_{4}  \partial_{4}$\\
 & $
   + x_{2} {{x_{3}}^2} \partial_{3}
   + x_{2}  x_{4} \partial_{5} + 2\, x_{3} x_{4}  \partial_{6}
   + 2\,  {{x_{2}}^2}  x_{3} x_{6} \partial_{4} + 2\,
x_{2} x_{3} x_{6}  \partial_{6} $\\
   \hline
   $V_3'$ & $ {{x_{1}}^2} x_{2}  \partial_{4}
    + 2\, {{x_{1}}^2} \partial_{6}
   + 2\,  x_{1} {{x_{2}}^2} \partial_{5}
    +x_{1}  x_{2}  x_{3}  \partial_{6}
     +  {{x_{2}}^2} {{x_{3}}^2}  \partial_{6} $\\
$V_3''$& $ {{x_{1}}^2}  x_{2}  \partial_{1}
 +{{x_{1}}^2}x_{5}  \partial_{4}
 +  {{x_{1}}^2}  \partial_{3}
 + x_{1}  x_{2}x_{3}  \partial_{3}
  + 2\,  x_{1}  x_{2}x_{5}  \partial_{5}
   +  x_{1} x_{3}  x_{5} \partial_{6}$\\
   &$+  x_{2}  {{x_{3}}^2} x_{5}  \partial_{6}
     +  x_{1} x_{2}  x_{4} \partial_{4}
      + 2\,  x_{1}  {{x_{2}}^2}  \partial_{2}
       + x_{1} x_{2}x_{3} \partial_{3}
        + 2\,  x_{1} x_{2} x_{6}  \partial_{6}$\\
        &$ + 2\, x_{1}  x_{4}  \partial_{6}
          + 2\, {{x_{2}}^2} {{x_{3}}^2}  \partial_{3}
           + 2\,  {{x_{2}}^2}x_{3}x_{6}  \partial_{6}
           + 2\, {{x_{2}}^2} x_{4} \partial_{5}
           +  x_{2}x_{3} x_{4}  \partial_{6} $\\
\hline $V_4'$ & $  {{x_{1}}^2}{{x_{2}}^2} \partial_{1}
        + {{x_{1}}^2} x_{2} x_{5} \partial_{4}
 + {{x_{1}}^2}  x_{2} \partial_{3}
 + 2\, {{x_{1}}^2} x_{5}\partial_{6}
 +  x_{1} {{x_{2}}^2} x_{3} \partial_{3}
 + 2\, x_{1}  {{x_{2}}^2} x_{5}  \partial_{5}$\\
 &$+ x_{1}  x_{2} x_{3}  x_{5} \partial_{6}
 + {{x_{2}}^2} {{x_{3}}^2}  x_{5} \partial_{6}   $\\
$V_4''$ & $ {{x_{1}}^2} x_{4}
     \partial_{4} +  {{x_{1}}^2}  x_{5} \partial_{5}
   +  {{x_{1}}^2}x_{6} \partial_{6} + {{x_{2}}^2} {{x_{3}}^2}x_{6} \partial_{6}
   +  x_{1} {{x_{2}}^2}{{x_{3}}^2} \partial_{1}
   + x_{1}  {{x_{2}}^2}x_{3}  \partial_{2}
   + 2\, x_{1}  {{x_{2}}^2} x_{6}  \partial_{5}$\\
   &$
   + 2\, x_{1}{{x_{3}}^2}  x_{5} \partial_{6}
   + 2\,  x_{1} x_{2} {{x_{3}}^2} \partial_{3}
   + 2\,  x_{1}  x_{2}  x_{4}\partial_{5} +  x_{1} x_{3} x_{4}
        \partial_{6} + x_{2} {{x_{3}}^2} x_{4} \partial_{6}
   + 2\, x_{1} {{x_{2}}^2}x_{3} x_{6} \partial_{4} $\\
              \hline
$V_5'$ & $ {{x_{1}}^2}  {{x_{2}}^2} \partial_{2}
   +  {{x_{1}}^2}x_{4} \partial_{6}
   + 2\,  {{x_{1}}^2} x_{2} x_{3}
          \partial_{3} + 2\,{{x_{1}}^2}  x_{2}  x_{4}
          \partial_{4} +  {{x_{1}}^2}  x_{2} x_{6}
        \partial_{6} + 2\, {{x_{2}}^2}  {{x_{3}}^2}  x_{4}
        \partial_{6}$\\
        &$
   + x_{1}  {{x_{2}}^2}  {{x_{3}}^2} \partial_{3}
   + x_{1}  {{x_{2}}^2} x_{4}  \partial_{5}
   + x_{1} {{x_{2}}^2}  x_{3} x_{6}\partial_{6} + 2\, x_{1}  x_{2} x_{3}  x_{4}
   \partial_{6}$\\
          \hline
$V_6'$ & ${{x_{1}}^2} {{x_{2}}^2} x_{4}
     \partial_{1} + {{x_{1}}^2}{{x_{2}}^2} x_{5}\partial_{2} + 2\,
{{x_{1}}^2}  {{x_{2}}^2} x_{6}  \partial_{3}
   + {{x_{1}}^2} x_{2} x_{4} \partial_{3} + 2\,  {{x_{1}}^2}x_{4}  x_{5}
          \partial_{6} $\\
          &$+ 2\,
 {{x_{1}}^2} x_{2}  x_{3}x_{5} \partial_{3} +  {{x_{1}}^2} x_{2} x_{4}x_{5} \partial_{4}
   + {{x_{1}}^2}x_{2} x_{5} x_{6} \partial_{6}
   +{{x_{2}}^2}{{x_{3}}^2}x_{4}  x_{5}\partial_{6}$\\
   &$ + x_{1}  {{x_{2}}^2} {{x_{3}}^2}x_{5}
        \partial_{3} + x_{1}{{x_{2}}^2}x_{3}  x_{4}\partial_{3} + 2\,  x_{1}
{{x_{2}}^2}  x_{4}  x_{5} \partial_{5} +x_{1}{{x_{2}}^2}x_{3}
x_{5}
 x_{6}  \partial_{6}+  x_{1} x_{2}  x_{3}  x_{4}  x_{5}  \partial_{6}$\\
\hline
\end{tabular}
$$\normalsize
For $N=(2,1,1)$, the lowest hight vectors are as in the table
above together with the following ones
$$ \footnotesize
\renewcommand{\arraystretch}{1.4}
\begin{tabular}{|l|l|} \hline
$\fg_i$ & lowest weight vectors\\
\hline
$V_4'''$ & $
{{x_{1}}^3}\partial_{4}
   + 2\, {{x_{1}}^2} x_{2}\partial{5}
   +  {{x_{1}}^2}  x_{3}\partial_{6}
   +  x_{1}  x_{2}{{x_{3}}^2} \partial_{6}  $\\
              \hline
$\dots$&$\dotfill$\\ \hline $V_{11}'$ & ${{x_{1}}^5} {{x_{2}}^2}
\partial_{2}
   +  {{x_{1}}^5} x_{4}  \partial_{6}
   +  {{x_{1}}^4}  {{x_{2}}^2}{{x_{3}}^2}  \partial_{3}
   + {{x_{1}}^4} {{x_{2}}^2}  x_{4}  \partial_{5}
   + 2\,  {{x_{1}}^5} x_{2}  x_{3} \partial_{3} + 2\,
 {{x_{1}}^5}  x_{2}  x_{4} \partial_{4}$\\
 &$ + {{x_{1}}^5} x_{2}  x_{6}
        \partial_{6} + 2\, {{x_{1}}^3}  {{x_{2}}^2}  {{x_{3}}^2}  x_{4} \partial_{6}
   + {{x_{1}}^4} {{x_{2}}^2}  x_{3}  x_{6} \partial_{6} + 2\,
 {{x_{1}}^4}  x_{2} x_{3}x_{4} \partial_{6}$\\
          \hline
$V_{12}'$ & $ {{x_{1}}^5}{{x_{2}}^2} x_{4} \partial_{1} +
{{x_{1}}^5} {{x_{2}}^2}  x_{5}  \partial_{2}
   + 2\, {{x_{1}}^5} {{x_{2}}^2}  x_{6}  \partial_{3}
   + {{x_{1}}^5}  x_{2} x_{4} \partial_{3} + 2\,
{{x_{1}}^5}  x_{4} x_{5} \partial_{6} +  {{x_{1}}^4} {{x_{2}}^2}
 {{x_{3}}^2}x_{5}  \partial_{3}$\\
&$   + {{x_{1}}^4} {{x_{2}}^2} x_{3}  x_{4}\partial_{3} + 2\,
{{x_{1}}^4}  {{x_{2}}^2} x_{4}  x_{5} \partial_{5} + 2\,
 {{x_{1}}^5} x_{2} x_{3} x_{5}\partial_{3} + {{x_{1}}^5} x_{2} x_{4} x_{5}
        \partial_{4} +{{x_{1}}^5}  x_{2}x_{5} x_{6}
        \partial_{6}$\\
        &$ +  {{x_{1}}^3}{{x_{2}}^2}{{x_{3}}^2} x_{4} x_{5}
        \partial_{6} +  {{x_{1}}^4}  {{x_{2}}^2}x_{3} x_{5} x_{6}
        \partial_{6} +  {{x_{1}}^4} x_{2} x_{3}x_{4} x_{5}  \partial_{6}$\\
\hline
\end{tabular}
$$\normalsize

Let us investigate if $\fBj(3;N|3)$ has partial prolongs as
subalgebras:

(i) Denote by $\fg_1'$ the $\fg_0$-module generated by $V_1'$. We
have $\sdim(\fg_1')=2|2$. The CTS partial prolong
$(\fg_{-},\fg_0,\fg_1')_*$ gives a graded Lie superalgebra with the
property that $[\fg_{-1},\fg_1]\simeq \{Y_1,h_1\}:=\mathfrak{aff}$.
From the description of irreducible modules over solvable Lie
superalgebras \cite{Ssol}, we see that the irreducible
$\mathfrak{aff}$-modules are 1-dimensional. For irreducible
$\mathfrak{aff}$-submodules  $\fg_{-1}'$ in $\fg_{-1}$ we have two
possibilities: to take $\fg_{-1}'=\{Y_4\}$ or $\fg_{-1}'=\{Y_7\}$;
for both of them, $\fg_{-1}'$ is purely odd and we can never get a
simple Cartan prolong.

(ii) Denote by $\fg_1''$ the $\fg_0$-module generated by $V_1''$.
We have $\sdim(\fg_1'')=2|2$. The CTS partial prolong
$(\fg_{-},\fg_0,\fg_1'')_*$ returns $\fosp(2|4)$.

\ssec{A description of $\fBj(2|4)$} We consider ${}^{3)}\fg(2|3)$
with $r=(1,0,0)$. In this case, $\sdim(\fg(2,3)_-)=2|4$. Since the
$\fg(2,3)_0$-module action is not faithful, we consider the quotient
algebra $\fg_0=\fg(2,3)_0/\mathfrak{ann}(\fg_{-1})$ and embed
$(\fg(2,3)_{-},\fg_0)\subset \fvect (2|4)$. This realization is
given by the following table:\tiny
 $$ \small
\renewcommand{\arraystretch}{1.4}
\begin{tabular}{|l|l|} \hline
$\fg_{i}$&the generators (even $\mid$ odd) \\ \hline \hline
$\fg_{-1}$&$Y_{6}=\partial_2,\; Y_{8}=\partial_1 \;\mid\;
Y_{11}=\partial_3, \; Y_{10}=\partial_4,\; Y_4=\partial_5, \;
Y_1=\partial_6\;$\\ \hline $\fg_{0}\simeq$& $Y_3=x_{2}\partial_{1} +
2\, x_{4}\partial_{3} + x_{6}\partial_{5}, Y_9=[Y_2,[Y_3,Y_5]],
Z_3=x_{1}\partial_{2} + 2\, x_{3}\partial_{4} + x_{5}\partial_{6},
Z_9=$\\ &$[Z_2, [Z_3,Z_5]],\; H_2=[Z_2, Y_2], H_3=[Z_3,Y_3] \;\mid\;
   Y_2=x_{1}\partial_{4} + x_{5}\partial_{2},Y_5=[Y_2,Y_3],$\\
$ \fosp(3|2)$& $Y_7=[Y_3, [Y_2,Y_3]],\; Z_2=
    x_{2}\partial_{5} + 2 x_{4}\partial_{1}, Z_5=[Z_2,Z_3],
   Z_7=[Z_3, [Z_2,Z_3]] $
   \\
\hline
\end{tabular}
$$\normalsize
The $\fg_0$-module $\fg_{-1}$ is irreducible, having one lowest
weight vector $Y_{11}$ and one highest weight vector $Y_1$. The CTS
prolong $(\fg_-,\fg_0)_*$ gives a Lie superalgebra of superdimension
$13|14$. Indeed, $\sdim(\fg_1)=4|4$ and $\sdim(\fg_2)=1|0$. The
$\fg_0$-module $\fg_1$ has one lowest vector: $$ \small
\renewcommand{\arraystretch}{1.4}
\begin{tabular}{|l|} \hline
$V_1=2\,x_{1}  x_{2} \partial_{3}
   + x_{1}x_{6}  \partial_{1} + 2\, {{x_{2}}^2}  \partial_{4}
   +  x_{2}  x_{5}  \partial_{1} + 2\,  x_{2}  x_{6}  \partial_{2}
   + x_{4}  x_{5} \partial_{3} + 2\,  x_{4}x_{6} \partial_{4}
   +  x_{5}x_{6}  \partial_{5}$\\
\hline
\end{tabular}
$$ The $\fg_2$ is one-dimensional spanned by the following vector $$
\small
\renewcommand{\arraystretch}{1.4}
\begin{tabular}{|l|} \hline $ 2\,  {{x_{1}}^2}  x_{2} \partial_{1}
    +  {{x_{1}}^2} x_{4}\partial_{3}
    + 2\,  {{x_{1}}^2} x_{6}  \partial_{5}
    +  x_{1}  {{x_{2}}^2}  \partial_{2}
    + 2\,  x_{1} x_{2} x_{3}  \partial_{3}
    +  x_{1}  x_{2}  x_{4}  \partial_{4}
    + 2\,  x_{1}  x_{2}  x_{5} \partial_{5}
    + x_{1} x_{2}  x_{6}  \partial_{6}$\\
    $    +  x_{1}  x_{3}  x_{6}  \partial_{1}
    + 2\,  x_{1}  x_{4}  x_{5}\partial_{1}
    +  x_{1} x_{4}  x_{6} \partial_{2}
    + 2\,{{x_{2}}^2} x_{3} \partial_{4}
    + {{x_{2}}^2}  x_{5} \partial_{6}
    +  x_{2} x_{3} x_{5}  \partial_{1}
    + 2\,  x_{2} x_{3}  x_{6} \partial_{2}
    +  x_{2}  x_{4}  x_{5} \partial_{2}$\\
    $ +  x_{3} x_{4}  x_{5} \partial_{3}
    + 2\, x_{3} x_{4} x_{6}  \partial_{4}
    + x_{3} x_{5}  x_{6}  \partial_{5}
    + 2\,  x_{4} x_{5} x_{6} \partial_{6}$\\
\hline
\end{tabular}
$$
Besides, if $i>2$, then $\fg_i=0$ for all values of the sharing
parameter $N=(N_1,N_2)$. A direct computation gives
$[\fg_1,\fg_1]=\fg_2$ and $[\fg_{-1},\fg_1]=\fg_0$. SuperLie tells
us that this Lie superalgebra has three ideals $I_1\subset
I_2\subset I_3$ with the same non-positive part but different
positive parts: $\sdim(I_1)=10|14$, $\sdim(I_2)=11|14$,
$\sdim(I_3)=12|14$. The ideal $I_1$ is just our $\fbj$, see
\cite{BjL, CE}. The partial CTS prolong with $I_1$ returns $I_1$
plus an  outer derivation given by the vector above (of degree 2).
It is clear now that $\fBj(2|4)$ is not simple.

\ssec{A description of $\fBj(3|5)$} We consider ${}^{3)}\fg(2|3)$
and $r=(0,1,0)$. In this case, $\sdim(\fg(2,3)_-)=3|5$. Since the
$\fg(2,3)_0$-module action is again not faithful, we consider the
quotient module $\fg_0=\fg(2,3)_0/\mathfrak{ann}(\fg_{-1})$ and
embed $(\fg(2,3)_{-},\fg_0)\subset \fvect (3;\underline{N}|5)$. This
realization is given by the following table:\tiny
$$ \small
\renewcommand{\arraystretch}{1.4}
\begin{tabular}{|l|l|} \hline
$\fg_{i}$&the generators (even $\mid$ odd) \\ \hline \hline
$\fg_{-2}$ & $Y_9=\partial_1\;\mid\; Y_{10}=\partial_3,
Y_{11}=\partial_2$ \\ \hline $\fg_{-1}$&$Y_{8}=\partial_4,\;
Y_{6}=\partial_5 \;\mid\; Y_5=2\, x_{4}\partial_{2} + 2\,
x_{5}\partial_{3} + 2\, x_{7}\partial_{1}
   + \partial_{7}, Y_2=x_{4}\partial_{3} - 2\, x_{6}\partial_{1} +
   \partial_{8}$\\
   &$ Y_7=x_{5}\partial_{2} + \partial_{6}$\\
\hline $\fg_{0}\simeq \fsl(1|2)$& $ H_1=[Z_1,Y_1],
H_3=[Z_3,Y_3],Y_3=2\, x_{3}\partial_{2} +
   2\,  x_{7}  x_{8} \partial_{1} + x_{5}\partial_{4}
   + 2\, x_{7}\partial_{6} + x_{8}\partial_{7},$\\
   & $ Z_3=2\, x_{2}\partial_{3} +
   2\, x_{6} x_{7} \partial_{1} + x_{4}\partial_{5}
   + x_{6}\partial_{7} + 2\, x_{7}\partial_{8}\;\mid\;
   Y_4=[Y_1,Y_3], \; Z_4=[Z_1,Z_3],$
   \\
   &$Y_1=2\, \left( 2\, x_{1}\partial_{3}
      + 2\, x_{6}  x_{7}  \partial_{2}
      + x_{6}\partial_{4} + x_{7}\partial_{5} \right)$\\
&$ Z_1=2\, \left( x_{3}\partial_{1} + 2\,
            x_{4}  x_{5}  \partial_{2}
      + 2\,  {{x_{5}}^2}  \partial_{3}
      + 2\,  x_{5}  x_{7}  \partial_{1}
      + 2\, x_{4}\partial_{6} + x_{5}\partial_{7} \right),
      $\\
\hline
\end{tabular}
$$\normalsize
The $\fg_0$-module $\fg_{-1}$ is irreducible, having one highest
weight vector $Y_2$. We have $\sdim(\fg_1)=6|4$. The $\fg_0$-module
$\fg_1$ has two lowest weight vectors given by
$$
\small
\renewcommand{\arraystretch}{1.4}
\begin{tabular}{|l|l|} \hline
$V_1'$ & $x_{1} x_{5}\partial_{2} + 2\, x_{5}
 x_{6} x_{8} \partial_{2} +  x_{5} x_{7} x_{8}  \partial_{3}
        + 2\, x_{1}\partial_{6} + 2\, x_{3}\partial_{4}
        +  x_{5}  x_{7} \partial_{4}
        + x_{5}  x_{8}  \partial_{5}
        + 2\, x_{7} x_{8}  \partial_{7}$\\
   $V_1''$ & $  x_{6}  x_{7}  x_{8}  \partial_{2}
           + 2\, x_{1}\partial_{4}
           +  x_{7}  x_{8}  \partial_{5} $\\
\hline
\end{tabular}
$$
Now, the $\fg_0$-module generated by the the vectors $V_1'$ and
$V_1''$ is not the whole $\fg_1$ but a $\fg_0$-module that we denote
by $\fg_1''$, of $\sdim=4|4$. The CTS prolong $(\fg_{-}, \fg_0,
\fg_1)_*$ is not simple, so consider the Lie subsuperalgebra
$(\fg_{-}, \fg_0, \fg_1'')_*$; the superdimensions of its positive
part are
$$
\renewcommand{\arraystretch}{1.4}
\begin{tabular}{|l|c|c|c|c|c|} \hline
$\mathrm{ad}_{\fg_1''}^i(\fg_1'')$ &$\fg_1''$ &
$\mathrm{ad}_{\fg_1''}(\fg_1'')$&
$\mathrm{ad}_{\fg_1''}^2(\fg_1'')$&$\mathrm{ad}_{\fg_1''}^3(\fg_1'')$
& $\mathrm{ad}_{\fg_1''}^4(\fg_1'')$ \\ \hline $\sdim$ &$4|4$ &
$4|4$ & $4|4$& $3|2$& $2|1$\\ \hline
\end{tabular}
$$
The lowest weight vectors of the above components are precisely
$\{V_2', V_2'', V_3, V_4, V_5\}$ described bellow: \tiny
$$ \small
\renewcommand{\arraystretch}{1.4}
\begin{tabular}{|l|l|} \hline
$\mathrm{ad}_{\fg_1}^i(\fg_1)$ & lowest weight vectors\\ \hline
$V_2'$ & $  {{x_{1}}^2}\partial_2
        + 2\,  x_{1}  x_{7} \partial_{4}
        + 2\, x_{1}  x_{8}  \partial_{5}
        +  x_{1}  x_{6}  x_{8}  \partial_{2}
        + 2\,  x_{1}  x_{7} x_{8}  \partial_{3} $\\
$V_2''$ &$2\, {{x_{1}}^2} \partial_{1}
           +  x_{1}  x_{2} \partial_{2}
           + x_{1} x_{3} \partial_{3}
           + x_{1}  x_{6}  \partial_{6}
           + x_{1} x_{7} \partial_{7}
           + x_{1}x_{8} \partial_{8}
           + 2\, x_{2}x_{7}  \partial_{4}
           + 2\, x_{2} x_{8} \partial_{5}$\\
           &$
           + 2\, x_{3} x_{6} \partial_{4}
           + 2\, x_{3} x_{7} \partial_{5}
           + x_{2}x_{6}x_{8} \partial_{2}
           + 2\, x_{2} x_{7}  x_{8} \partial_{3}
           + x_{3} x_{6} x_{7}  \partial_{2}
           + x_{6}  x_{7} x_{8}  \partial_{7} $\\
           \hline
           $V_3$ & ${{x_{1}}^2}\partial_{4}
   + 2\,  x_{1}x_{7}x_{8}\partial_{5}
   + 2\,  x_{1}  x_{6}  x_{7}  x_{8} \partial_{2}$\\
   \hline
   $V_4$ & $ {{x_{1}}^2}  x_{3} \partial_{2}
   + 2\,  {{x_{1}}^2}  x_{5} \partial_{4}
   + {{x_{1}}^2}  x_{7} \partial_{6}
   + 2\, {{x_{1}}^2}  x_{8} \partial_{7}
   +  {{x_{1}}^2}  x_{7}  x_{8} \partial_{1}
   + 2\, x_{1} x_{3}  x_{7}\partial_{4}
   + 2\,  x_{1} x_{3} x_{8} \partial_{5}$\\
   &$
   + x_{1} x_{3}x_{6}  x_{8}\partial_{2}
   + 2\, x_{1}x_{3}x_{7} x_{8}\partial_{3}
   + x_{1} x_{5}  x_{7} x_{8} \partial_{5}
   + x_{1}  x_{5} x_{6}  x_{7}x_{8}\partial_{2}$\\
\hline $V_5$ & ${{x_{1}}^2}  x_{2} \partial_{4}
   + 2\, {{x_{1}}^2} x_{3} \partial_{5}
   + 2\, {{x_{1}}^2}  x_{6}  x_{7}\partial_{6}
   +  {{x_{1}}^2} x_{6} x_{8} \partial_{7}
   + 2\, {{x_{1}}^2} x_{7}  x_{8}\partial_{8}
   + 2\, {{x_{1}}^2} x_{6}x_{7} x_{8}
   \partial_{1}$\\
   &$   + 2\,  x_{1} x_{2}  x_{7} x_{8}\partial_{5}
   + 2\, x_{1} x_{3} x_{6} x_{7}\partial_{4}
   + 2\, x_{1} x_{3} x_{6} x_{8} \partial_{5}
   + 2\, x_{1}  x_{2}x_{6} x_{7} x_{8}\partial_{2}
   + 2\,x_{1}x_{3} x_{6} x_{7}  x_{8} \partial_{3}$\\
   \hline
\end{tabular}
$$\normalsize
Since none of the known simple finite dimensional Lie superalgebra
over (algebraically closed) fields of characteristic 0 or $>3$ has
such a non-positive part in any $\Zee$-grading, it follows that
$\fBj(3; \underline{N}|5)$ is new.

Let us investigate if $\fBj(3;\underline{N}|5)$ has subalgebras
--- partial prolongs.

(i) Denote by $\fg_1'$ the $\fg_0$-module generated by $V_1'$. We
have $\sdim(\fg_1')=2|3$. The CTS partial prolong
$(\fg_{-1},\fg_0,\fg_1')_*$ gives a graded Lie superalgebra with
$\sdim (\fg_2')=2|2$ and $\fg_i'=0$ for $i>3$. An easy computation
shows that $[\fg_{-1}, \fg_1']=\fg_0$ and $[\fg_1',\fg_1']\subsetneq
\fg_2'$. Since we are investigating simple Lie superalgebra, we take
the simple part of $(\fg_{-1},\fg_0,\fg_1')_*$. This simple Lie
superalgebra is isomorphic to $\fg(2, 3)/\fc=\fbj$.

(ii) Denote by $\fg_1''$ the $\fg_0$-module generated by $V_1''$. We
just saw that $\sdim(\fg_1'')=4|4$. The CTS partial prolong
$(\fg_{-1},\fg_0,\fg_1'')_*$ gives also $\fBj(3|5)$.

\underline{$r=(0,0,1)$}. In this case, $\sdim(\fg(2,3)_-)=4|5$.
Since the $\fg(2,3)_0$-module action is not faithful, we consider
the quotient algebra $\fg_0=\fg(2,3)_0/\mathfrak{ann}(\fg_{-1})$ and
embed $(\fg(2,3)_{-},\fg_0)\subset \fvect (4;\underline{N}|5)$. The
CTS prolong returns $\fbj:=\fg(2,3)/\fc$.

\ssec{A description of $\mathfrak{Me}(3; N|3)$}

1) Our first idea was to try to repeat the above construction with a
suitable super version of $\fg(2)$. There is only one simple super
analog of $\fg(2)$, namely $\fag(2)$, but our attempts \cite{BjL} to
construct a super analog of Melikyan algebra in the above way as
Kuznetsov suggested \cite{Ku1} (reproduced in \cite{GL4}) resulted
in something quite distinct from the Melikyan algebra: The Lie
superalgebras we obtained, an exceptional one $\fbj$ (cf. \cite{CE,
BGL1}) and a series $\fBj$, are indeed simple but do not resemble
either  $\fg(2)$ or $\mathfrak{Me}$.

2) Our other idea is based on the following observation. The
anti-symmetric form
\begin{equation}\label{pairing}
(f, g):=\int fdg =\int fg' dt,
\end{equation} on the quotient space $F/\text{const}$ of
functions (with compact support) modulo constants on the
1-dimensional manifolds, has its counterpart in $1|1$-dimensional
case in presence of a contact structure \so{and only in this case}
as follows from the description of invariant bilinear differential
operators, see \cite{KLV}. Indeed, the Lie superalgebra $\fk(1|1)$
does not distinguish between the space of volume forms (let its
generator be denoted $\vvol$) and the quotient $\Omega^1/F\alpha$,
where $\alpha=dt+\theta d\theta$ is the contact form.

For any prime $p$ therefore, on the space $\fg_{-1}:=\cO(1;
\underline{N}| 1)/\const$ of \lq\lq functions (with compact support)
in one even indeterminate $u$ and one odd, $\theta$  modulo
constants", the superanti-symmetric bilinear form
\begin{equation}\label{spairing}
(f, g):=\int (fdg \mod F\alpha)=\int (f_0g_0'-f_1g_1) dt,
\end{equation}
where $f=f_0(t)+f_1(t)\theta$ and $g=g_0(t)+g_1(t)\theta$ and
where $':=\frac{d}{dt}$, is nondegenerate.

Therefore, we may expect that, for $p$ small and
$\underline{N}=1$, the Melikyan effect will reappear. Consider
$p=5$ as the most plausible.

We should be careful with parities. The parity of $\vvol$ is a
matter of agreement, let it be even. Then the integral is an odd
functional but the factorization modulo $F\alpha$ makes the form
$(\ref{spairing})$ even. (Setting $p(\vvol)=\od$ we make the
integral an even functional and the factorization modulo $F\alpha$
makes the form $(\ref{spairing})$ even again.)

Since the form $(\ref{spairing})$ is even, we get the following
realization of $$\fk(1;\underline{1}|1)\subset
\fosp(5|4)\simeq\fk(5;\underline{1, ..., 1}|5)$$ by generating
functions of contact vector fields on the $5|5$-dimensional
superspace with the contact form, where the coefficients are found
from the explicit values of
$$
d\hat t+ \sum_ i(\hat p_i d\hat q_i-\hat q_2 d\hat p_i) + \sum_j
\left (\hat \xi_j d\eta_j +\hat \eta_j d\hat \xi_j\right) - \hat
\theta d\hat \theta.
$$
The coordinates on this $5|5$-dimensional superspace are hatted in
order not to confuse them with generating functions of
$\fk(1;\underline{1}|1)$:
\begin{equation}\label{smeme}
\renewcommand{\arraystretch}{1.4}
\begin{tabular}{|c|l|}
\hline $\fg_{i}$&basis elements\\
\hline $\fg_{-2}$&$\hat 1$\\
\hline $\fg_{-1}$&$\hat p_1=t,\; \hat p_2=t^2,\;  \hat q_1=t^3,\;
\hat q_2=t^4,$\\
&$\hat \xi_1=\theta,\;   \hat \xi_2=t\theta,\;  \hat
\theta=t^2\theta,\;
\hat \eta_2=t^3\theta, \; \hat \eta_1=t^4\theta$\\
\hline \end{tabular} \end{equation} We explicitly have:
\begin{equation}\label{explsuper}
\renewcommand{\arraystretch}{1.4}
\begin{array}{l}
(t, t^4)=\int_{\underline{N}} t\cdot
t^3dt_{\underline{N}}=\int_{\underline{N}}4t^4dt_{\underline{N}}=4=-(t^4, t);\\
(t^2, t^3)=\int_{\underline{N}} t^2\cdot
t^2dt_{\underline{N}}=\int_{\underline{N}}6t^4dt_{\underline{N}}=1=-(t^3, t^2);\\
(t^4\theta, \theta)=-\int_{\underline{N}} t^4\cdot 1dt_{\underline{N}}=-1=(\theta, t^4\theta);\\
(t^3\theta, t\theta)=-\int_{\underline{N}} t^3\cdot tdt_{\underline{N}}=-4=(t\theta, t^3\theta);\\
(t^2\theta, t^2\theta)=-\int_{\underline{N}} t^2\cdot
t^2dt_{\underline{N}}=-6=-1.
\end{array}
\end{equation}
Now, let us realize $\fk(1;\underline{1}|1)$ by contact fields in
hatted functions:
\begin{equation}\label{smeme1}
\renewcommand{\arraystretch}{1.4}
\begin{tabular}{|c|l|}
\hline $\fg_{i}$&basis elements\\
\hline $\fg_{-2}$&$\hat 1$\\
\hline $\fg_{-1}$&$\hat p_1=t,\; \hat p_2=t^2,\;  \hat q_2=4
t^3,\;
\hat q_1=t^4,$\\
&$\xi_1=\theta,\;   \xi_2=t\theta,\;  \hat \theta=t^2\theta,\;
\eta_2=4 t^3\theta, \; \eta_1=t^4\theta$\\
\hline $\fg_{0}$&$1=2\,  \hat p_{1}  {\cdot}
      \hat q_{2}   + 2{{\hat p_{2}}^2}
   + 3\,  \xi_1 \eta_2+ 3\,  \xi_2 \hat \theta ;\quad
   t=2\,  \hat p_{1} \hat q_{1} + 4\,  \hat p_{2} \hat q_{2}
   + 4\,  \xi_1  \eta_1+ 2\,  \xi_2\eta_2; $\\
&$t^2=2\,  \hat p_{2} \hat q_{1}
   + 4\, {{\hat q_{2}}^2}
   + 4\, \xi_2  \eta_1
   + \hat \theta  \eta_2; \quad t^3= 3\, \hat q_{1}  \hat q_{2}
   + 4\, \hat \theta\eta_1; \quad t^4={{\hat q_{1}}^2} +  \eta_2  \eta_1 ;$\\
       & $ \theta= \hat p_{1}  \eta_2
   +  \hat p_{2} \hat \theta
   +  \hat q_{1} \xi_1
   + \hat q_{2}  \xi_2 ;\quad  t\theta= \hat p_{1}  \eta_1
   + 2\,  \hat p_{2}  \eta_2
   + \hat q_{1} \xi_2
   + 2\, \hat q_{2}  \hat \theta ;$\\
         & $ t^2\theta= \hat p_{2} \eta_1
   +  \hat q_{1} \hat \theta + 2\,  \hat q_{2}  \eta_2; \quad
   t^3\theta= 4\, \hat q_{1}\eta_2
   + 4\, \hat q_{2} \eta_1 ;\quad  t^4 \theta= \hat q_{1}  \eta_1 $\\
\hline \end{tabular}
\end{equation}
The CTS prolong gives that $\fg_1=0$.

The case where $p=3$ is more interesting because it will give us the
series $\mathfrak{Me}(3;N|3)$. The non-positive part is as follows:
\begin{equation}
\renewcommand{\arraystretch}{1.4}
\begin{tabular}{|c|l|}
\hline $\fg_{i}$&basis elements\\
\hline $\fg_{-2}$&$\hat 1$\\
\hline $\fg_{-1}$&$\hat p_1=t,\; \hat q_2=t^2,\;\xi_1=\theta,\;
\hat \theta=t\theta,\;
\eta_1=t^2\theta$\\
\hline $\fg_{0}$&$1=\hat p_1^2+2\hat \xi_1 \hat \eta_1 ;\quad
   t=2\, \hat p_{1}  \hat q_{1}
   + 2\,  \hat \xi_{1}  \hat \eta_{1};\,
   t^2=2\, {{\hat q_{1}}^2} + 2\,  \hat \theta
        \hat \eta_{1} ; \; \theta= 2\, \hat p_1 \hat \theta
   + \hat q_1  \hat \xi_1 ;$\\
&$    t\theta=  \hat p_1\hat \eta_1
   +  \hat q_1  \hat \theta ;\; t^2\theta=  \hat q_1  \hat \eta_1 $\\
\hline \end{tabular}
\end{equation}

The Lie superalgebra $\fg_0$ is not simple because
$[\fg_{-1},\fg_1]=\fg_0\backslash \{ t^2\theta=  \hat q_1  \hat
\eta_1 \}$. Denote $\fg_0':=[\fg_{-1},\fg_1]\simeq \fosp(1|2)$. The
CTS partial prolong $(\fg_{-},\fg_0')_*$ seems to be very
interesting. First, our computation shows that the parameter
$\underline{M}=(M_1,M_2,M_3)$ depends only on the first parameter
(relative to $t$). Namely, $\underline{M}=(M_1,1,1)$. For $M_1=1$,
2, the super dimensions of the positive components of $\fBj(3;M|3)$
are given in the following table:
$$
\renewcommand{\arraystretch}{1.4}
\begin{tabular}{|l|c|c|c|c|c|c|c|c|c|c|} \hline
$M=(1,1,1)$ &$\fg_1'$ & $\fg_2'$& $\fg_3'$&$\fg_4'$ & $\fg_5'$ &
--&--&--&--&--\\
\hline $\sdim$ &$2|4$ & $4|2$& $2|4$& $3|2$& $0|1$& &&&&\\
\hline
 \hline
$M=(2,1,1)$ &$\fg_1'$ & $\fg_2'$& $\fg_3'$&$\fg_4'$ & $\fg_5'$
&$\cdots$ & $\fg_{14}'$ & $\fg_{15}'$& $\fg_{16}'$ & $\fg_{17}'$\\
\hline $\sdim$ &$2|4$
& $4|2$& $2|4$& $4|2$ & $2|4$ &$\cdots$& $4|2$ & $2|4$ &$3|2$ & $0|1$\\
\hline
\end{tabular}
$$
Here we have that $[\fg_{-1},\fg_1]=\fg_0'$ and the $\fg_1'$
generates the positive part. The standard criteria for simplicity
ensures that $\mathfrak{Me}(3;N|3)$ is simple. For $N=(1,a,b)$, the
lowest weight vectors are as follows:\footnotesize
$$
\renewcommand{\arraystretch}{1.4}
\begin{tabular}{|l|l|} \hline
$\fg_i$ & lowest weight vectors\\ \hline $V_1'$ & $ 2{{p_{1}}^{(2)}}
\eta_{1}
        + 2\, p_{1} q_{1}  \theta
        +  {{q_{1}}^{(2)}}  \xi_{1}
        + \xi_{1}  \theta  \eta_{1}   $\\
$V_1''$ &$  2\, p_{1}  q_{1} \eta_{1}
           +   {{q_{1}}^{(2)}}\theta +   \eta_{1}t  $\\
   \hline
   $V_2'$ & $ 2\,  p_{1}  q_{1}  \theta \eta_{1}
        +  {{q_{1}}^{(2)}}\xi_{1} \eta_{1}
        +  t  {{q_{1}}^{(2)}}
        +  t  \theta  \eta_{1}  $\\
$V_2''$& $  2\,
              {{p_{1}}^{(2)}} {{q_{1}}^{(2)}}
           +  {{p_{1}}^{(2)}} \theta \eta_{1}
           + p_{1} q_{1}  \xi_{1}  \eta_{1}
           + 2 {{q_{1}}^{(2)}}  \xi_{1} \theta
           + {t^{(2)}}$\\
\hline $V_3'$ & $  2\,  t  {{p_{1}}^{(2)}}
 \eta_{1}
        +  t  p_{1}  q_{1}  \theta
        + 2\,  t {{q_{1}}^{(2)}}  \xi_{1}
        +  t \xi_{1}  \theta \eta_{1}$\\
$V_3''$ & $  2\, {{p_{1}}^{(2)}} {{q_{1}}^{(2)}}  \eta_{1}
           + 2\, t p_{1}  q_{1}  \eta_{1}
           + 2\, {{q_{1}}^{(2)}}  \xi_{1}\theta  \eta_{1}
           +  t  {{q_{1}}^{(2)}} \theta
           + {t^{(2)}} \eta_{1}$\\
              \hline
$V_4'$ & $   2\,  {{p_{1}}^{(2)}}
           {{q_{1}}^{(2)}}\theta  \eta_{1}
      + 2\,  t  p_{1} q_{1}  \theta  \eta_{1}
      +  t  {{q_{1}}^{(2)}}\xi_{1}  \eta_{1}
      +  {t^{(2)}}{{q_{1}}^{(2)}}
      +  {t^{(2)}} \theta  \eta_{1} $\\
          \hline
$V_5'$ & $  {{p_{1}}^{(2)}}  {{q_{1}}^{(2)}} \xi_{1}  \theta
\eta_{1}
      +   {t^{(2)}}  {{p_{1}}^{(2)}} \eta_{1}
      +  {t^{(2)}}  p_{1} q_{1}  \theta
      + 2\, {t^{(2)}} {{q_{1}}^{(2)}}  \xi_{1}
      + 2\,  {t^{(2)}} \xi_{1} \theta  \eta_{1} $\\
\hline
\end{tabular}
$$\normalsize
For $N=(2,a,b)$, the lowest weight vectors are as above together
with:
$$ \footnotesize
\renewcommand{\arraystretch}{1.4}
\begin{tabular}{|l|l|} \hline
$\fg_i$ & lowest weight vectors\\ \hline $V_4''$ & $  2\,
            t {{p_{1}}^{(2)}}  {{q_{1}}^{(2)}}
                +  t {{p_{1}}^{(2)}}
               \theta  \eta_{1}
           +  t p_{1} q_{1} \xi_{1}  \eta_{1}
           + 2\, t  {{q_{1}}^{(2)}}
               \xi_{1}  \theta
           + {t^{(3)}}    $\\
   \hline
   $V_5''$ & $  2\,t  {{p_{1}}^{(2)}} {{q_{1}}^{(2)}} \eta_{1}
           + 2\,  {t^{(2)}} p_{1}  q_{1}  \eta_{1}
           + 2\, t  {{q_{1}}^{(2)}}\xi_{1} \theta \eta_{1}
           +  {t^{(2)}} {{q_{1}}^{(2)}}
               \theta
           + {t^{(3)}} \eta_{1} $\\\hline
          $\dots$&$\dotfill$\\ \hline
$V_{15}'$& $2\,{t^{(5)}} {{p_{1}}^{(2)}} {{q_{1}}^{(2)}}
            \xi_{1} \theta\eta_{1}
        + 2{t^{(7)}} {{p_{1}}^{(2)}} \eta_{1}
        + 2\,{t^{(7)}}  p_{1}  q_{1}\theta
        + {t^{(7)}}{{q_{1}}^{(2)}} \xi_{1}
        + {t^{(7)}}  \xi_{1} \theta \eta_{1}  $\\
& $  2\, {t^{(6)}}  {{p_{1}}^{(2)}} {{q_{1}}^{(2)}} \eta_{1}
           + 2\, {t^{(7)}} p_{1}q_{1} \eta_{1}
           + 2\, {t^{(6)}} {{q_{1}}^{(2)}} \xi_{1}  \theta \eta_{1}
           +   {t^{(7)}}{{q_{1}}^{(2)}} \theta
           + {t^{(8)}} \eta_{1}   $\\
\hline $V_{16}''$& $2\, {t^{(6)}} {{p_{1}}^{(2)}} {{q_{1}}^{(2)}}
\theta \eta_{1}
      + 2\,  {t^{(7)}}  p_{1}  q_{1}  \theta\eta_{1}
      +  {t^{(7)}}{{q_{1}}^{(2)}}  \xi_{1} \eta_{1}
      + {t^{(8)}}{{q_{1}}^{(2)}}
      + {t^{(8)}}  \theta  \eta_{1}   $\\
\hline $V_{17}''$ & $ 2 {t^{(6)}} {{p_{1}}^{(2)}} {{q_{1}}^{(2)}}
\xi_{1}  \theta  \eta_{1}
      + 2\,  {t^{(8)}}{{p_{1}}^{(2)}}  \eta_{1}
      +  2\, {t^{(8)}}p_{1} q_{1} \theta
      +   {t^{(8)}} {{q_{1}}^{(2)}}  \xi_{1}
      + {t^{(8)}}  \xi_{1} \theta  \eta_{1} $\\
\hline
\end{tabular}
$$
Let us investigate the subalgebras of $\mathfrak{Me}(3;N|3)$ ---
partial prolongs:

(i) Denote by $\fg_1'$ the $\fg_0$-module generated by $V_1'$. We
have $\sdim(\fg_1')=0|1$ and $\fg_i=0$ for all $i>1$. The CTS
partial prolong $(\fg_{-},\fg_0,\fg_1')_*$ gives a graded Lie
superalgebra with the property that $[\fg_{-1},\fg_1]\simeq
\fosp(1|2)$. The partial CTS prolong $(\fg_{-},\fosp(1|2))_*$ is
not simple

(ii) Denote by $\fg_1''$ the $\fg_0$-module generated by $V_1''$. We
have $\sdim(\fg_1'')=3|2$. The CTS partial prolong
$(\fg_{-},\fg_0,\fg_1'')_*$ returns $\mathfrak{brj}(2;3)$.

\ssec{A description of $\mathfrak{Brj}(4|3)$} We have the following
realization of the non-positive part inside $\fvect(4|3)$:
\begin{equation}\label{Brj(4|3)}
\tiny
\renewcommand{\arraystretch}{1.4}
\begin{tabular}{|l|l|} \hline
$\fg_{i}$&the generators (even $\mid$ odd) \\ \hline \hline
$\fg_{-4}$&$Y_{8}=\partial_1\;\mid\; Y_{7}=\partial_5$\\ \hline
$\fg_{-3}$&$Y_{6}=\partial_2\;\mid $\\ \hline
$\fg_{-2}$&$Y_{4}=\partial_3\;\mid\; Y_{5}=x_{3}\partial_{5} +
x_{6}\partial_{1} + \partial_{6}$\\ \hline $\fg_{-1}$&$Y_3=2\,
x_{2}\partial_{1} + 2\, x_{3}\partial_{2} + \partial_{4} \;\mid$\\
&$ Y_2=x_{2}\partial_{5} + 2\, {{x_{4}}^{(2)}} x_{7}
\partial_{1}
   + x_{4} x_{6} \partial_{1} +  x_{6} x_{7} \partial_{5}
   + x_{4}  x_{7} \partial_{2} + x_{6}\partial_{2} + 2\, x_{4}\partial_{6}
   + 2\, x_{7}\partial_{3} + 1\partial_{7},$\\ \hline $\fg_{0}\simeq \mathfrak{hei}(0|2)\oplus \mathbb{K}$&
$H_1=[Z_1,Y_1],\;H_2=2\, x_{5}\partial_{5} + x_{2}\partial_{2}
   + 2\, x_{3}\partial_{3} + 2\, x_{4}\partial_{4}
   + x_{7}\partial_{7},\;\mid$\\
   & $ Y_1=2\,  {{x_{3}}^{(2)}}  \partial_{5}
   + 2\, x_{3}  x_{6} \partial_{1} + 2\, x_{5}\partial_{1}
   + 2\, x_{3}\partial_{6} + x_{7}\partial_{4},$\\
   &$ Z_1= {{x_{4}}^{(2)}} \partial_{6}
   + 2\, {{x_{4}}^{(2)}} x_{6}  \partial_{1}
   +  {{x_{4}}^{(2)}}x_{7} \partial_{2}
   + x_{4}  x_{7} \partial_{3} + 2\,  x_{4} x_{6}  x_{7} \partial_{5} + x_{1}\partial_{5}
   + 2\, x_{4}\partial_{7} + 2\, x_{6}\partial_{3} $
   \\
\hline
\end{tabular}
\end{equation}
The Lie superalgebra $\fg_0$ is solvable, and hence the CTS prolong
$(\fg_{-},\fg_0)_*$ is NOT simple since $\fg_1$ does not generate
the positive part. Our calculation shows that the prolong does not
depend on $N$, i.e., $N=(1,1,1,1)$. The simple part of this prolong
is $\mathfrak{brj}(2;3)$. The $\sdim$ of the positive parts are
described as follows:
$$
\renewcommand{\arraystretch}{1.4}
\begin{tabular}{|l|c|c|c|c|c|c|c|c|c|c|} \hline
&$\fg_1$ & $\fg_2$& $\fg_3$&$\fg_4$ & $\fg_5$ & $\fg_6$&$\fg_7$&$\fg_8$&$\fg_9$
&$\fg_{10}$\\
\hline $\sdim$ &$1|1$ & $2|2$& $1|2$& $2|2$& $1|1$ & $2|2$ & $1|1$ & $1|1$ &
 $0|1$&$1|1$\\
\hline
\end{tabular}
$$
and the lowest weight vectors are as follows:
$$ \tiny
\renewcommand{\arraystretch}{1.4}
\begin{tabular}{|l|l|} \hline
$\fg_i$ & lowest weight vectors\\ \hline $V_1'$ & $  2\, x_{2}
x_{3} \partial_{5}
    + 2\,  x_{2} x_{6} \partial_{1}
    +  x_{3}  x_{4} \partial_{6}    + 2\,  x_{3}  x_{6} \partial_{2}
    +  x_{3}  x_{7} \partial_{3}    +  x_{4} x_{7}\partial_{4}
    + x_{3} {{x_{4}}^{(2)}} x_{7}\partial_{1}
    + 2\, x_{3} x_{4} x_{6}\partial_{1}$ \\
    & $   + 2\,  x_{3}  x_{4} x_{7} \partial_{2}
    + 2\,  x_{3}x_{6}x_{7}\partial_{5} + 2\, x_{2}\partial_{6}
    + 2\, x_{3}\partial_{7}    + 2\, x_{5}\partial_{2} + 2\, x_{6}\partial_{4}$\\
   \hline
   $V_2'$ & $ {{x_{4}}^{(2)}}\partial_{4}
    +  x_{1}  x_{3} \partial_{5}    +  x_{1} x_{6} \partial_{1}
    + x_{3} {{x_{4}}^{(2)}} \partial_{6} + 2\,
        x_{3}  x_{4} \partial_{7} + 2\,  x_{3} x_{6}\partial_{3}
    + 2\,  x_{4}  x_{6} \partial_{4}
    + 2\, x_{6} x_{7} \partial_{7}$\\
    &$ + 2\,  x_{3} {{x_{4}}^{(2)}} x_{6} \partial_{1} +
    x_{3}  {{x_{4}}^{(2)}} x_{7} \partial_{2}
    +  x_{3}x_{4}x_{7}\partial_{3}+ 2\, x_{3} x_{4}  x_{6}  x_{7} \partial_{5}
    + x_{1}\partial_{6} + 2\, x_{5}\partial_{3} $\\
    $V_2''$
    &
    $   2\,  {{x_{2}}^{(2)}} \partial_{1} +  {{x_{3}}^{(2)}} \partial_{3}
    + 2\, x_{2} x_{3} \partial_{2}
    + x_{3}  x_{4} \partial_{4}    +  x_{3}  x_{5}\partial_{5}
    + 2\, x_{3} x_{7}\partial_{7}    + x_{5}  x_{6} \partial_{1}
       + 2\,x_{6}  x_{7} \partial_{4} + x_{2}\partial_{4}
       + x_{5}\partial_{6}$\\\hline
$V_{3}'$& $2\,  x_{1}  x_{2} \partial_{1}
    + 2\, x_{1}  x_{3} \partial_{2}
    +  x_{2} x_{3} \partial_{3}    + 2\,  x_{2}  x_{4} \partial_{4}
    + 2\, x_{2}  x_{5}\partial_{5}
    + x_{2}  x_{6} \partial_{6}    + 2\,  x_{2}x_{7} \partial_{7}
    + 2\,  x_{3} x_{6}\partial_{7}    +  x_{4}  x_{5}
    \partial_{6}$\\
&$  + 2\, x_{5}  x_{6} \partial_{2}    +  x_{5}  x_{7}
\partial_{3}
    +  {{x_{4}}^{(2)}}  x_{5} x_{7}\partial_{1}
    + x_{3}  x_{4}  x_{6}\partial_{6}    + x_{3} x_{4} x_{7} \partial_{7}
    + x_{3}  x_{6}  x_{7} \partial_{3}    + 2\, x_{4}x_{5}  x_{6} \partial_{1}
    + 2\, x_{4}  x_{5}x_{7}\partial_{2}$\\
    &$  + 2\,  x_{4} x_{6} x_{7} \partial_{4}
    + 2\, x_{5}  x_{6}  x_{7}\partial_{5}
     +  x_{3}  {{x_{4}}^{(2)}}x_{6} x_{7}\partial_{1}
    + 2\, x_{3}  x_{4} x_{6} x_{7}\partial_{2}  + x_{1}\partial_{4}
    + 2\, x_{5}\partial_{7} $\\
$V_{3}''$& $ 2\, {{x_{3}}^{(2)}} \partial_{7} + {{x_{3}}^{(2)}}
 x_{4}\partial_{6} + 2\,{{x_{3}}^{(2)}} \partial x_{6} \partial_{2}
    + {{x_{3}}^{(2)}}  x_{7} \partial_{3}
    + 2\,  x_{2} {{x_{3}}^{(2)}} \partial_{5}
    + 2\,  x_{2}  x_{3} \partial_{6}
    + 2\, x_{2}  x_{5} \partial_{1}    + x_{2}  x_{7}
    \partial_{4}$\\
    &$
    + 2\,  x_{3} x_{5} \partial_{2}    + 2\, x_{3} x_{6} \partial_{4}
    +  {{x_{3}}^{(2)}}  {{x_{4}}^{(2)}}  x_{7} \partial_{1}
    + 2\, {{x_{3}}^{(2)}} x_{4}  x_{6} \partial_{1}
    + 2\, {{x_{3}}^{(2)}} x_{4} x_{7} \partial_{2}
    + 2\,{{x_{3}}^{(2)}} x_{6} x_{7}\partial_{5} + 2\,x_{2} x_{3}
    x_{6}\partial_{1}$\\
    &$ + x_{3}  x_{4} x_{7}\partial_{4}  + x_{5}\partial_{4}$\\
\hline
$V_{4}$& $ 2\, {{x_{2}}^{(2)}} \partial_{6} +
2\,{{x_{2}}^{(2)}}  x_{3} \partial_{5}
    + 2\, {{x_{2}}^{(2)}}  x_{6} \partial_{1}
      + 2\, {{x_{3}}^{(2)}} {{x_{4}}^{(2)}} \partial_{6}
    +  {{x_{3}}^{(2)}}  x_{4} \partial_{7}
    + {{x_{3}}^{(2)}} x_{6}\partial_{3}
    + 2\, x_{1}  {{x_{3}}^{(2)}}\partial_{5}$\\
    &$ + 2\,  x_{1}  x_{3} \partial_{6}
    + 2\,  x_{1} x_{5} \partial_{1}
    +  x_{1}  x_{7} \partial_{4}
    + 2\,  x_{2} x_{3} \partial_{7}
    + 2\,  x_{2}  x_{5} \partial_{2}
    + 2\, x_{2} x_{6} \partial_{4}
    + 2\, x_{3}  x_{5} \partial_{3}
    +  x_{5}  x_{6} \partial_{6}
    +  x_{5}  x_{7} \partial_{7}$\\
    &$ + {{x_{3}}^{(2)}}  {{x_{4}}^{(2)}}  x_{6} \partial_{1}
    + 2\, {{x_{3}}^{(2)}}  {{x_{4}}^{(2)}} x_{7} \partial_{2}
    + 2\, {{x_{3}}^{(2)}} x_{4}  x_{7} \partial_{3}
    + 2\, x_{1} x_{3} x_{6} \partial_{1}
    +  x_{2}  x_{3}x_{4} \partial_{6}
    + 2\, x_{2}x_{3} x_{6} \partial_{2}$\\
    &$+ x_{2}  x_{3} x_{7}\partial_{3}
    +  x_{2}  x_{4}x_{7}\partial_{4}
    + x_{3}  {{x_{4}}^{(2)}} x_{7}\partial_{4}
    + 2\,  x_{3} x_{4}  x_{6} \partial_{4}
    +  {{x_{3}}^{(2)}}  x_{4} x_{6}  x_{7} \partial_{5}
    + x_{2}  x_{3}{{x_{4}}^{(2)}} x_{7} \partial_{1}
    + 2\, x_{2} x_{3}  x_{4} x_{6} \partial_{1}$\\
    &$+ 2\,  x_{2}  x_{3} x_{4}  x_{7} \partial_{2}
    + 2\,  x_{2} x_{3}  x_{6}  x_{7} \partial_{5} $\\
           $V_4''$ & $  {{x_{1}}^{(2)}} \partial_{1}
    +  {{x_{4}}^{(2)}} x_{5}\partial_{6}
    + 2\, x_{1}  x_{3} \partial_{3}
    +  x_{1}  x_{4} \partial_{4}
    + x_{1}  x_{5} \partial_{5}
    + 2\, x_{1}  x_{6} \partial_{6}
    + x_{1}  x_{7} \partial_{7}
    + 2\,  x_{4} x_{5} \partial_{7}
    + 2\, x_{5} x_{6} \partial_{3}$\\
    &$  + 2\,  {{x_{4}}^{(2)}}x_{5}  x_{6} \partial_{1}
    +  {{x_{4}}^{(2)}} x_{5} x_{7}\partial_{2}
    + 2\,  {{x_{4}}^{(2)}}  x_{6}  x_{7}\partial_{4}
    + x_{3}  {{x_{4}}^{(2)}}  x_{6}\partial_{6}
    + x_{3}  {{x_{4}}^{(2)}} x_{7}\partial_{7}
    + 2\, x_{3}  x_{4}x_{6} \partial_{7}$\\
    &$ +  x_{4}  x_{5}x_{7} \partial_{3}
    + x_{3}  {{x_{4}}^{(2)}}  x_{6} x_{7}\partial_{2}
    + x_{3}  x_{4} x_{6}x_{7} \partial_{3}
    + 2\, x_{4}  x_{5} x_{6}  x_{7}\partial_{5} $\\
\hline $V_{5}'$ & $ x_{1} x_{2} \partial_{6}
    +  x_{1}  x_{3}\partial_{7}
    +  x_{1} x_{5}\partial_{2}
    + x_{1} x_{6} \partial_{4}
    + 2\,  x_{2}  x_{5} \partial_{3}
    +  x_{5}  x_{6} \partial_{7}
    + x_{1} x_{2}  x_{3}\partial_{5}
    + x_{1}  x_{2}x_{6}\partial_{1}
    + 2\, x_{1}  x_{3} x_{4}\partial_{6}$\\
    &$ + x_{1}  x_{3}x_{6}\partial_{2}
    + 2\, x_{1}  x_{3}  x_{7}\partial_{3}
    + 2\, x_{1} x_{4} x_{7}\partial_{4}
    +  x_{2}  {{x_{4}}^{(2)}} x_{7}\partial_{4}
    +  x_{2}  x_{3} {{x_{4}}^{(2)}}\partial_{6} + 2\,
          x_{2}  x_{3} x_{4}\partial_{7}
    + 2\,  x_{2}  x_{3} x_{6} \partial_{3}$\\
&$ + 2\,  x_{2}  x_{4} x_{6}\partial_{4}
    + 2\,  x_{2} x_{6}  x_{7} \partial_{7}
    + 2\, x_{4}  x_{5} x_{6} \partial_{6}
    + 2\, x_{4}  x_{5} x_{7} \partial_{7}
    + 2\,  x_{5} x_{6} x_{7} \partial_{3}
    + 2\, {{x_{4}}^{(2)}} x_{5}  x_{6}  x_{7}\partial_{1}$\\
&$ + 2\, x_{1}  x_{3}  {{x_{4}}^{(2)}} x_{7}\partial_{1}
    + x_{1}  x_{3}  x_{4} x_{6} \partial_{1}
    + x_{1} x_{3}  x_{4} x_{7}\partial_{2}
    + x_{1}  x_{3}  x_{6} x_{7} \partial_{5}
    + 2\,  x_{2}  x_{3} {{x_{4}}^{(2)}}  x_{6} \partial_{1}
    +  x_{2}  x_{3}  {{x_{4}}^{(2)}} x_{7} \partial_{2}$\\
&$ + x_{2}  x_{3}  x_{4} x_{7}\partial_{3}
    + 2\, x_{3}  x_{4} x_{6}  x_{7}\partial_{7}
    + x_{4}  x_{5}  x_{6}  x_{7}\partial_{2}
    + 2\, x_{2}  x_{3}x_{4}  x_{6} x_{7}\partial_{5} $\\
\hline $V_6'$ & ${{x_{1}}^{(2)}} \partial_{6} + {{x_{1}}^{(2)}}
x_{3}\partial_{5} +  {{x_{1}}^{(2)}}  x_{6}
\partial_{1} + 2\, x_{1}  x_{5} \partial_{3}
    +  {{x_{4}}^{(2)}}  x_{5} x_{6} \partial_{6}
    +  {{x_{4}}^{(2)}}  x_{5} x_{7} \partial_{7}
    + x_{1}  {{x_{4}}^{(2)}}x_{7} \partial_{4}$\\
&$ + x_{1} x_{3}{{x_{4}}^{(2)}} \partial_{6} + 2\,
     x_{1}  x_{3} x_{4} \partial_{7}
    + 2\, x_{1}  x_{3} x_{6} \partial_{3}
    + 2\,x_{1}  x_{4} x_{6}\partial_{4}
    + 2\,  x_{1}  x_{6}  x_{7}\partial_{7}
    + 2\,  x_{4}  x_{5}  x_{6} \partial_{7}$\\
&$    +  {{x_{4}}^{(2)}}  x_{5} x_{6} x_{7} \partial_{2}
    + 2\, x_{1}  x_{3}  {{x_{4}}^{(2)}}  x_{6}\partial_{1}
    + x_{1}x_{3} {{x_{4}}^{(2)}} x_{7} \partial_{2}
    +  x_{1} x_{3} x_{4}x_{7}\partial_{3}
    +  x_{3}  {{x_{4}}^{(2)}}  x_{6} x_{7}\partial_{7}$\\
    &$    + x_{4} x_{5} x_{6}x_{7}\partial_{3}
    + 2\,  x_{1}  x_{3} x_{4} x_{6}  x_{7}\partial_{5}$\\
$V_6''$ & ${{x_{2}}^{(2)}}  x_{3}\partial_{3}
    + 2\, {{x_{2}}^{(2)}}  x_{4} \partial_{4}
      + 2\,{{x_{2}}^{(2)}} x_{5} \partial_{5}
    + {{x_{2}}^{(2)}} x_{6} \partial_{6}
    + 2\, {{x_{2}}^{(2)}} x_{7} \partial_{7}
      + 2\,  x_{1} {{x_{2}}^{(2)}} \partial_{1}
    + x_{1}  {{x_{3}}^{(2)}}\partial_{3}$\\
    &$+ x_{1} x_{2}\partial_{4}
    + x_{1} x_{5} \partial_{6}
    + 2\, x_{2}  x_{5} \partial_{7}
    + 2\, {{x_{3}}^{(2)}} {{x_{4}}^{(2)}}x_{6}\partial_{6}
    + 2\, {{x_{3}}^{(2)}}  {{x_{4}}^{(2)}} x_{7}\partial_{7}
    +  {{x_{3}}^{(2)}} x_{4}x_{6} \partial_{7}
    + 2\, x_{1}  x_{2} x_{3} \partial_{2}$\\
&$    +  x_{1}  x_{3} x_{4} \partial_{4}
    +  x_{1}  x_{3} x_{5} \partial_{5}
    + 2\, x_{1} x_{3} x_{7} \partial_{7}
    +  x_{1} x_{5}x_{6} \partial_{1}
    + 2\,  x_{1} x_{6}  x_{7} \partial_{4}
    + 2\,  x_{2} x_{3} x_{6} \partial_{7}
    + x_{2}  x_{4} x_{5} \partial_{6}$\\
&$ + 2\, x_{2}  x_{5} x_{6} \partial_{2}
    +  x_{2}  x_{5}x_{7}\partial_{3}
    + x_{3}  {{x_{4}}^{(2)}}x_{5}\partial_{6}
    + 2\, x_{3} x_{4}  x_{5} \partial_{7}
    + 2\, x_{3} x_{5}  x_{6}\partial_{3}
    + x_{5}  x_{6}  x_{7}\partial_{7}$\\
&$  + 2\,  {{x_{3}}^{(2)}}  {{x_{4}}^{(2)}}  x_{6}
x_{7}\partial_{2}
    + 2\,  {{x_{3}}^{(2)}}  x_{4} x_{6} x_{7}\partial_{3}
    + x_{2} {{x_{4}}^{(2)}} x_{5}  x_{7}\partial_{1}
    + x_{2}  x_{3} x_{4}x_{6} \partial_{6}
    + x_{2}  x_{3}  x_{4} x_{7}\partial_{7}
    + x_{2}  x_{3}  x_{6} x_{7} \partial_{3}$\\
&$ + 2\,  x_{2} x_{4} x_{5}  x_{6}\partial_{1}
    + 2\, x_{2}  x_{4}x_{5}x_{7} \partial_{2}
    + 2\,  x_{2}  x_{4}x_{6} x_{7} \partial_{4}
    + 2\, x_{2}  x_{5} x_{6} x_{7} \partial_{5}
    + 2\, x_{3} {{x_{4}}^{(2)}} x_{5} x_{6} \partial_{1}
    +  x_{3} {{x_{4}}^{(2)}}x_{5}x_{7} \partial_{2}$\\
&$    + 2\,  x_{3}{{x_{4}}^{(2)}} x_{6} x_{7}\partial_{4}
    + x_{3}  x_{4}  x_{5}x_{7} \partial_{3}
    +x_{2}  x_{3} {{x_{4}}^{(2)}} x_{6} x_{7} \partial_{1}
    + 2\,  x_{2}  x_{3}  x_{4} x_{6} x_{7}\partial_{2}
    + 2\, x_{3} x_{4} x_{5} x_{6} x_{7} \partial_{5} $\\
\hline
\end{tabular}
$$

$$ \tiny
\renewcommand{\arraystretch}{1.4}
\begin{tabular}{|l|l|}
\hline $V_7'$& ${{x_{1}}^{(2)}}\partial_{4} + 2\, {{x_{1}}^{(2)}}
 x_{2} \partial_{1}
    + 2\, {{x_{1}}^{(2)}} x_{3} \partial_{2}
      + 2\, x_{1}  x_{5} \partial_{7}
    + x_{1} x_{2}  x_{3} \partial_{3}
    + 2\,  x_{1} x_{2} x_{4} \partial_{4}
    + 2\, x_{1} x_{2} x_{5} \partial_{5}
    +  x_{1}  x_{2} x_{6} \partial_{6}$\\
&$    + 2\, x_{1}  x_{2}x_{7} \partial_{7}
    + 2\,  x_{1}  x_{3}  x_{6} \partial_{7}
    +  x_{1} x_{4}  x_{5} \partial_{6}
    + 2\, x_{1} x_{5} x_{6} \partial_{2}
    + x_{1}  x_{5}x_{7} \partial_{3}
    + 2\,  x_{2}  {{x_{4}}^{(2)}}  x_{5} \partial_{6}
    + x_{2}  x_{4} x_{5} \partial_{7}$\\
&$ +x_{2}  x_{5} x_{6}\partial_{3}
    + x_{1}  {{x_{4}}^{(2)}}  x_{5} x_{7}\partial_{1}
    + x_{1}  x_{3}x_{4}  x_{6} \partial_{6}
    + x_{1}  x_{3} x_{4} x_{7}\partial_{7}
    + x_{1}  x_{3} x_{6} x_{7}\partial_{3}
    + 2\,x_{1}  x_{4} x_{5} x_{6} \partial_{1}
    + 2\, x_{1} x_{4}  x_{5}x_{7} \partial_{2}$\\
&$    + 2\,  x_{1}  x_{4} x_{6}  x_{7}\partial_{4}
    + 2\, x_{1}x_{5}x_{6} x_{7} \partial_{5}
    +  x_{2}  {{x_{4}}^{(2)}} x_{5} x_{6}\partial_{1}
    + 2\, x_{2}  {{x_{4}}^{(2)}}  x_{5}  x_{7} \partial_{2}
    +  x_{2}  {{x_{4}}^{(2)}} x_{6}  x_{7} \partial_{4}
    + 2\,  x_{2}  x_{3}  {{x_{4}}^{(2)}}  x_{6} \partial_{6}$\\
    &$+ 2\,
        x_{2} x_{3}  {{x_{4}}^{(2)}} x_{7} \partial_{7} +
        x_{2}  x_{3}x_{4}x_{6}\partial_{7}
    + 2\, x_{2}  x_{4} x_{5} x_{7} \partial_{3}
    + x_{4}  x_{5} x_{6} x_{7} \partial_{7}
    + x_{1}  x_{3}  {{x_{4}}^{(2)}} x_{6} x_{7} \partial_{1}
    + 2\, x_{1}  x_{3}  x_{4} x_{6}  x_{7} \partial_{2}$\\
&$    + 2\, x_{2} x_{3} {{x_{4}}^{(2)}}  x_{6} x_{7} \partial_{2}
    + 2\, x_{2}  x_{3} x_{4} x_{6}  x_{7}\partial_{3}
    +  x_{2} x_{4}  x_{5} x_{6}  x_{7} \partial_{5}$\\
$V_8'$ & ${{x_{1}}^{(2)}} {{x_{3}}^{(2)}}\partial_{5} +
 {{x_{1}}^{(2)}}  x_{3} \partial_{6} +  {{x_{1}}^{(2)}}  x_{5}\partial_{1} + 2\,
         {{x_{1}}^{(2)}}  x_{7} \partial_{4}
    + 2\,  {{x_{2}}^{(2)}} x_{5} \partial_{3}
      + x_{1}  {{x_{2}}^{(2)}}\partial_{6} +
       {{x_{1}}^{(2)}}x_{3}  x_{6}\partial_{1}$\\
&$    + {{x_{2}}^{(2)}} {{x_{4}}^{(2)}} x_{7} \partial_{4}
    + {{x_{2}}^{(2)}}  x_{3} {{x_{4}}^{(2)}} \partial_{6} + 2\,
         {{x_{2}}^{(2)}} x_{3} x_{4} \partial_{7}
    + 2\, {{x_{2}}^{(2)}}  x_{3}  x_{6} \partial_{3}
    + 2\, {{x_{2}}^{(2)}}  x_{4}  x_{6}\partial_{4}
    + 2\, {{x_{2}}^{(2)}}  x_{6}x_{7} \partial_{7}$\\
&$    + x_{1}  {{x_{2}}^{(2)}}  x_{3}\partial_{5}
    + x_{1} {{x_{2}}^{(2)}}x_{6}\partial_{1}
    + x_{1} {{x_{3}}^{(2)}} {{x_{4}}^{(2)}}\partial_{6} + 2\,
         x_{1}  {{x_{3}}^{(2)}}  x_{4}\partial_{7} + 2\,
          x_{1}  {{x_{3}}^{(2)}} x_{6} \partial_{3} +
           x_{1}  x_{2}  x_{3}\partial_{7}
    + x_{1}  x_{2}x_{5}\partial_{2}$\\
&$    + x_{1}  x_{2} x_{6} \partial_{4}
    +  x_{1}  x_{3}  x_{5}\partial_{3}
    + 2\, x_{1}  x_{5}  x_{6} \partial_{6}
    + 2\,  x_{1} x_{5} x_{7}\partial_{7}
    + x_{2}  x_{5}  x_{6}\partial_{7}
    + 2\,  {{x_{2}}^{(2)}} x_{3} {{x_{4}}^{(2)}}x_{6} \partial_{1} +
     {{x_{2}}^{(2)}}x_{3}{{x_{4}}^{(2)}} x_{7}\partial_{2}$\\
&$    + {{x_{2}}^{(2)}} x_{3}  x_{4}  x_{7}\partial_{3}
    +{{x_{3}}^{(2)}} {{x_{4}}^{(2)}}x_{6} x_{7} \partial_{7}
    + 2\, x_{1} {{x_{3}}^{(2)}}  {{x_{4}}^{(2)}}  x_{6} \partial_{1}
    +x_{1}  {{x_{3}}^{(2)}}  {{x_{4}}^{(2)}}  x_{7} \partial_{2}
    +  x_{1}  {{x_{3}}^{(2)}}x_{4}  x_{7} \partial_{3}$\\
&$ + 2\,  x_{1}  x_{2}  x_{3} x_{4} \partial_{6}
    +  x_{1}  x_{2}  x_{3}  x_{6}\partial_{2}
    + 2\, x_{1} x_{2}x_{3} x_{7}\partial_{3}
    + 2\,  x_{1} x_{2}  x_{4}  x_{7} \partial_{4}
    + 2\,  x_{1}  x_{3}  {{x_{4}}^{(2)}}  x_{7} \partial_{4}
    + x_{1}  x_{3}  x_{4}  x_{6} \partial_{4}$\\
&$    + 2\, x_{2}  x_{4}  x_{5} x_{6}\partial_{6}
    + 2\,  x_{2} x_{4} x_{5} x_{7} \partial_{7}
    + 2\, x_{2}  x_{5}  x_{6}  x_{7}\partial_{3}
    + 2\, x_{3}  {{x_{4}}^{(2)}} x_{5}  x_{6} \partial_{6}
    + 2\, x_{3}  {{x_{4}}^{(2)}}  x_{5}  x_{7}\partial_{7}
    + x_{3}  x_{4}  x_{5}x_{6} \partial_{7}$\\
&$    + 2\,  {{x_{2}}^{(2)}}x_{3} x_{4}  x_{6}  x_{7} \partial_{5}
    + 2\, x_{1}  {{x_{3}}^{(2)}}  x_{4} x_{6} x_{7}\partial_{5}
    + 2\, x_{1}  x_{2} x_{3} {{x_{4}}^{(2)}} x_{7} \partial_{1}
    + x_{1}  x_{2} x_{3} x_{4}x_{6} \partial_{1}
    +  x_{1}  x_{2} x_{3} x_{4}  x_{7}\partial_{2}$\\
&$    +  x_{1} x_{2}x_{3}x_{6}x_{7} \partial_{5}
    + 2\, x_{2}  {{x_{4}}^{(2)}}  x_{5}  x_{6} x_{7}\partial_{1}
    + 2\, x_{2} x_{3}  x_{4}  x_{6}  x_{7} \partial_{7}
    + x_{2}  x_{4}  x_{5} x_{6}  x_{7} \partial_{2}
    + 2\,  x_{3}  {{x_{4}}^{(2)}}  x_{5}  x_{6}  x_{7} \partial_{2}
    +$\\&$ 2\, x_{3}  x_{4}  x_{5} x_{6}  x_{7}\partial_{3}$\\\hline
    $V_9'$ & $  {{x_{1}}^{(2)}} x_{2}x_{3} \partial_{5}
    +{{x_{1}}^{(2)}}  x_{2}  x_{6}\partial_{1}
    + 2\,  {{x_{1}}^{(2)}}x_{3}  {{x_{4}}^{(2)}}  x_{7} \partial_{1} +
     {{x_{1}}^{(2)}} x_{3}  x_{4} x_{6}\partial_{1}
    + {{x_{1}}^{(2)}}  x_{3}x_{6}  x_{7} \partial_{5}
    + 2\,  x_{1}x_{2}  x_{3} {{x_{4}}^{(2)}}x_{6}\partial_{1}$\\
    &$+ 2\,
          x_{1}  x_{2} x_{3} x_{4} x_{6} x_{7} \partial_{5}
    + 2\,  x_{1} {{x_{4}}^{(2)}} x_{5}  x_{6}  x_{7} \partial_{1}
    + {{x_{1}}^{(2)}} x_{3} x_{4}  x_{7} \partial_{2}
    +  {{x_{1}}^{(2)}}  x_{3} x_{6} \partial_{2}
    +{{x_{1}}^{(2)}} x_{5} \partial_{2}$\\
&$  + x_{1}  x_{2}  x_{3} {{x_{4}}^{(2)}} x_{7}\partial_{2}
    +  x_{1} x_{4}  x_{5} x_{6}  x_{7} \partial_{2}
    + x_{2}  {{x_{4}}^{(2)}} x_{5}x_{6} x_{7} \partial_{2}
    + {{x_{1}}^{(2)}} x_{2} \partial_{6}
    + 2\,  {{x_{1}}^{(2)}}  x_{3} x_{4}\partial_{6}
    + 2\, {{x_{1}}^{(2)}}x_{3} x_{7} \partial_{3}$\\
&$ +  x_{1}  x_{2}  x_{3}{{x_{4}}^{(2)}} \partial_{6} +
     x_{1}  x_{2}x_{3}  x_{4}  x_{7}\partial_{3}
    + 2\,x_{1}  x_{2}  x_{3} x_{6} \partial_{3}
    + 2\, x_{1} x_{2}x_{5} \partial_{3}
    + 2\, x_{1}  x_{4}  x_{5} x_{6} \partial_{6}
    + 2\,  x_{1}  x_{5}x_{6}  x_{7}\partial_{3}$\\
&$ +  x_{2}  {{x_{4}}^{(2)}} x_{5}x_{6}\partial_{6}
    +  x_{2}  x_{4}  x_{5}x_{6} x_{7}\partial_{3}
    + {{x_{1}}^{(2)}} x_{3}\partial_{7}
    + 2\,  {{x_{1}}^{(2)}} x_{4} x_{7} \partial_{4}
    +  {{x_{1}}^{(2)}}  x_{6}\partial_{4}
    + 2\, x_{1}  x_{2}  x_{3} x_{4} \partial_{7}$\\
&$  +  x_{1} x_{2} {{x_{4}}^{(2)}} x_{7} \partial_{4}
    + 2\,  x_{1}  x_{2}  x_{4}x_{6}\partial_{4}
    + 2\, x_{1}  x_{2}x_{6}  x_{7} \partial_{7}
    + 2\, x_{1} x_{3} x_{4} x_{6} x_{7} \partial_{7}
    + 2\, x_{1}  x_{4} x_{5} x_{7} \partial_{7}
    +  x_{1}  x_{5}  x_{6} \partial_{7}$\\
&$ + x_{2}  x_{3}  {{x_{4}}^{(2)}}  x_{6}  x_{7} \partial_{7}
    + x_{2}  {{x_{4}}^{(2)}} x_{5} x_{7} \partial_{7}
    + 2\, x_{2} x_{4} x_{5} x_{6}\partial_{7} $\\
    \hline
    $V_{10}'$ & $  {{x_{1}}^{(2)}}  {{x_{2}}^{(2)}} \partial_{1}
    + 2\,  {{x_{1}}^{(2)}}  {{x_{3}}^{(2)}} \partial_{3}
    + 2\,  {{x_{1}}^{(2)}}  x_{2} \partial_{4}
      + 2\, {{x_{1}}^{(2)}}x_{5} \partial_{6}
    + {{x_{1}}^{(2)}} x_{2} x_{3} \partial_{2}
    + 2\,  {{x_{1}}^{(2)}}  x_{3}x_{4}\partial_{4}$\\
&$+ 2\, {{x_{1}}^{(2)}}  x_{3} x_{5}\partial_{5}
    + {{x_{1}}^{(2)}} x_{3} x_{7}\partial_{7}
    + 2\, {{x_{1}}^{(2)}} x_{5} x_{6}\partial_{1}
    +  {{x_{1}}^{(2)}}  x_{6}  x_{7}\partial_{4}
    + {{x_{2}}^{(2)}}  {{x_{4}}^{(2)}} x_{5}\partial_{6}
    + 2\,  {{x_{2}}^{(2)}} x_{4}x_{5} \partial_{7}$\\
    &$ + 2\,  {{x_{2}}^{(2)}}x_{5} x_{6} \partial_{3}
    + 2\,  x_{1}  {{x_{2}}^{(2)}} x_{3}\partial_{3} +
    x_{1}{{x_{2}}^{(2)}} x_{4} \partial_{4}
    +  x_{1}  {{x_{2}}^{(2)}} x_{5} \partial_{5}
    + 2\,  x_{1}  {{x_{2}}^{(2)}} x_{6} \partial_{6} +
   x_{1} {{x_{2}}^{(2)}} x_{7} \partial_{7}$\\
   &$ + x_{1} x_{2}x_{5} \partial_{7}
    + 2\, {{x_{2}}^{(2)}} {{x_{4}}^{(2)}}  x_{5}  x_{6}\partial_{1}
    + {{x_{2}}^{(2)}}  {{x_{4}}^{(2)}}x_{5}x_{7} \partial_{2}
    + 2\,  {{x_{2}}^{(2)}} {{x_{4}}^{(2)}}  x_{6}x_{7} \partial_{4}
    + {{x_{2}}^{(2)}}  x_{3}  {{x_{4}}^{(2)}}x_{6}\partial_{6}$\\
&$+ {{x_{2}}^{(2)}} x_{3}  {{x_{4}}^{(2)}} x_{7}\partial_{7}
    + 2\, {{x_{2}}^{(2)}}x_{3}x_{4} x_{6}\partial_{7}
    + {{x_{2}}^{(2)}} x_{4} x_{5}  x_{7}\partial_{3}
    + x_{1}  {{x_{3}}^{(2)}}  {{x_{4}}^{(2)}} x_{6}\partial_{6}
    + x_{1}  {{x_{3}}^{(2)}} {{x_{4}}^{(2)}}x_{7}  \partial_{7}$\\
&$ + 2\, x_{1}  {{x_{3}}^{(2)}} x_{4}  x_{6}\partial_{7}
    +x_{1}  x_{2}  x_{3} x_{6} \partial_{7}
    + 2\, x_{1}x_{2} x_{4} x_{5} \partial_{6}
    +  x_{1} x_{2} x_{5}  x_{6} \partial_{2}
    + 2\, x_{1} x_{2}  x_{5} x_{7}\partial_{3}
    + 2\, x_{1}  x_{3}  {{x_{4}}^{(2)}}  x_{5} \partial_{6}$\\
    &$+ x_{1} x_{3}  x_{4}x_{5} \partial_{7}
    + x_{1}  x_{3} x_{5} x_{6}\partial_{3}
    + 2\,  x_{1} x_{5} x_{6}  x_{7} \partial_{7}
    +{{x_{2}}^{(2)}} x_{3} {{x_{4}}^{(2)}} x_{6}  x_{7} \partial_{2}
    + {{x_{2}}^{(2)}} x_{3}  x_{4} x_{6}  x_{7} \partial_{3}$\\
&$+ 2\,  {{x_{2}}^{(2)}}  x_{4}  x_{5} x_{6} x_{7}
\partial_{5}
    + x_{1}  {{x_{3}}^{(2)}} {{x_{4}}^{(2)}} x_{6}x_{7}\partial_{2}
    +  x_{1} {{x_{3}}^{(2)}} x_{4}  x_{6}x_{7} \partial_{3}
    + 2\,  x_{1} x_{2} {{x_{4}}^{(2)}} x_{5}  x_{7}\partial_{1}
    + 2\,  x_{1}  x_{2} x_{3} x_{4}  x_{6} \partial_{6}$\\
&$    + 2\, x_{1} x_{2}  x_{3}  x_{4} x_{7} \partial_{7}
    + 2\, x_{1} x_{2} x_{3}  x_{6} x_{7} \partial_{3}
    +  x_{1} x_{2} x_{4} x_{5} x_{6} \partial_{1}
    + x_{1}  x_{2}  x_{4} x_{5} x_{7}\partial_{2}
    + x_{1} x_{2}  x_{4} x_{6}x_{7} \partial_{4}
    + x_{1}  x_{2}x_{5} x_{6}  x_{7} \partial_{5}$\\
&$ +  x_{1}  x_{3} {{x_{4}}^{(2)}}  x_{5}  x_{6} \partial_{1}
    + 2\, x_{1}  x_{3} {{x_{4}}^{(2)}}  x_{5} x_{7}\partial_{2}
    +x_{1}  x_{3}  {{x_{4}}^{(2)}} x_{6}  x_{7} \partial_{4}
    + 2\, x_{1} x_{3} x_{4}x_{5}  x_{7}\partial_{3}
    + 2\, x_{2}  x_{4} x_{5}  x_{6}  x_{7} \partial_{7}$\\
&$+ 2\,  x_{3} {{x_{4}}^{(2)}}x_{5} x_{6} x_{7} \partial_{7}
    + 2\,  x_{1}  x_{2}  x_{3} {{x_{4}}^{(2)}} x_{6} x_{7}\partial_{1}
    + x_{1}  x_{2} x_{3} x_{4}x_{6}x_{7} \partial_{2}
    +  x_{1}  x_{3}  x_{4} x_{5}  x_{6}  x_{7} \partial_{5}$\\\hline
\end{tabular}
$$
\normalsize

\ssec{A description of $\mathfrak{Brj}(3; \underline{N}|4)$} We have
the following realization of the non-positive part inside $\fvect(3;
\underline{N}|4)$:
\begin{equation}\label{Brj(3|4)}
\small
\renewcommand{\arraystretch}{1.4}
\begin{tabular}{|l|l|} \hline
$\fg_{i}$&the generators (even $\mid$ odd) \\ \hline \hline
$\fg_{-3}$&$\mid Y_{6}=\partial_4 $\\ \hline
$\fg_{-2}$&$Y_{5}=\partial_1,\; Y_{6}=\partial_2,\; Y_7=\partial_3\;\mid$\\
\hline $\fg_{-1}$&$\mid Y_2=2\, x_{3}\partial_{4} +
\partial_{5},\; Y_3= x_{2}\partial_{4} + x_{6}\partial_{1}
   + \partial_{6}$\\
&$Y_4=2\, x_{1}\partial_{4}
   + 2\, x_{5} x_{7} \partial_{4} + x_{5}\partial_{1}
   + x_{6}\partial_{2}+ 2\, x_{7}\partial_{3}+ \partial_{7}$\\
   \hline $\fg_{0}\simeq \mathfrak{hei}(2|0)\subplus \Kee H_2$& $H_1=[Z_1,Y_1],\;H_2=2\, x_{1}\partial_{1}
   + x_{3}\partial_{3}
   + x_{4}\partial_{4}
   + x_{6}\partial_{6}
   + 2\, x_{7}\partial_{7}$\\
   & $ Y_1=2\, x_{5} x_{6} x_{7} \partial_{4}
     + 2\, x_{1}\partial_{2}
   + 2\, x_{2}\partial_{3}   + 2\,  x_{5} x_{6}\partial_{1} +x_{6}x_{7}\partial_{3}
   + 2\, x_{5}\partial_{6}
   + 2\, x_{6}\partial_{7},$\\
   &$ Z_1= 2\,  x_{2}\partial_{1}
      + 2\, x_{3}\partial_{2}
      + x_{6}  x_{7} \partial_{1} + 2\, x_{6}\partial_{5}
      + x_{7}\partial_{6}\;\mid $
   \\
\hline
\end{tabular}
\end{equation}
The Lie superalgebra $\fg_0$ is solvable with the property that
$[\fg_0,\fg_0]=\mathfrak{hei}(2|0)$. The CTS prolong
$(\fg_{-},\fg_0)_*$ is NOT simple since $\fg_1$ does not generate
the positive part. Our calculation shows that the prolong does not
depend on $N$, i.e., $N=(1,1,1,1)$. The simple part of this prolong
is $\mathfrak{brj}$. The $\sdim$ of the positive parts are described
as follows:
$$
\renewcommand{\arraystretch}{1.4}
\begin{tabular}{|l|c|c|c|} \hline
&$\fg_1$ & $\fg_2$& $\fg_3$\\
\hline $\sdim$ &$0|3$ & $3|0$& $0|2$\\
\hline
\end{tabular}
$$
and the lowest weight vectors are
$$ \tiny
\renewcommand{\arraystretch}{1.4}
\begin{tabular}{|l|l|}
\hline $V_1'$& $  2\,
        {{x_{1}}^{(2)}} \partial_{4}
      + 2\,  x_{1}  x_{5} x_{7}\partial_{4}
      + x_{1}  x_{5} \partial_{1} +  x_{1}  x_{6}\partial_{2}
      + 2\,  x_{1} x_{7} \partial_{3} + x_{4}\partial_{3}
      + x_{1}\partial_{7}
      + 2\,  x_{5}  x_{6}\partial_{6}
      + x_{5} x_{7}\partial_{7}$\\\hline
$V_2'$ & $  2\, x_{1} x_{4} \partial_{4} +  x_{2}  x_{5}x_{6}
 x_{7} \partial_{4} + 2\,x_{4} x_{5} x_{7}\partial_{4}
      + 2\, {{x_{1}}^{(2)}} \partial_{1}
      + x_{1}  x_{2}\partial_{2} +  {{x_{2}}^
             {(2)}} \partial_{3}
        + x_{2} x_{5}  x_{6} \partial_{1} + 2\,
           x_{2}  x_{6} x_{7}\partial_{3}$\\
&$ + x_{4}  x_{5} \partial_{1} + x_{4}x_{6}\partial_{2}
      + 2\, x_{4}  x_{7} \partial_{3}
      + 2\,  x_{1}  x_{5}\partial_{5}
      +  x_{1}  x_{6} \partial_{6} + x_{2}  x_{5}\partial_{6}
      +x_{2} x_{6} \partial_{7} + x_{4}\partial_{7}
      + x_{5}  x_{6}  x_{7}\partial_{6}  $\\\hline
    $V_3'$ & $  {{x_{1}}^{(2)}}  x_{2} \partial_{4}
    + x_{1}  x_{2} x_{5} x_{7} \partial_{4} + 2\,
        x_{4}  x_{5} x_{6} x_{7} \partial_{4} +  {{x_{1}}^
             {(2)}}  x_{6}\partial_{1} + 2\,
         x_{1} x_{2}  x_{5} \partial_{1} + 2\,
         x_{1}  x_{2} x_{6}\partial_{2} +  x_{1} x_{2} x_{7}\partial_{3}$\\
         &$ + 2\,
          x_{1} x_{4} \partial_{2} + 2\,
          x_{1} x_{5} x_{6}x_{7}\partial_{1} + 2\,
          x_{2}  x_{4} \partial_{3} + 2\,
         x_{4}  x_{5} x_{6}\partial_{1} + x_{4}  x_{6} x_{7}\partial_{3}
          +  {{x_{1}}^
            {(2)}} \partial_{6}
    + 2\,  x_{1}  x_{2} \partial_{7}$\\
    &$ +  x_{1} x_{5}  x_{6}\partial_{5} + 2\,
         x_{1}  x_{5}x_{7} \partial_{6} +  x_{2}  x_{5}  x_{6} \partial_{6} + 2\,
       x_{2}  x_{5}  x_{7}\partial_{7} + 2\,
          x_{4} x_{5}\partial_{6} + 2\,
          x_{4}  x_{6} \partial_{7}$\\
    $V_{3}''$ & $ {{x_{1}}^
           {(2)}}  x_{3} \partial_{4}
    + x_{1}  {{x_{2}}^{(2)}}\partial_{4}
    + x_{1} x_{3}  x_{5} x_{7}\partial_{4} + {{x_{2}}^
             {(2)}}x_{5}  x_{7} \partial_{4} +  x_{1}x_{2} x_{6} \partial_{1} + 2\,
        x_{1}  x_{3}  x_{5} \partial_{1} + 2\,
         x_{1}  x_{3} x_{6} \partial_{2}$\\
         &$ +x_{1}  x_{3} x_{7}\partial_{3} + 2\,
         x_{1} x_{4} \partial_{1} + 2\,
          {{x_{2}}^{(2)}}  x_{5} \partial_{1}
    + 2\, {{x_{2}}^{(2)}} x_{6} \partial_{2}
     + {{x_{2}}^{(2)}}  x_{7} \partial_{3}
    + 2\, x_{2} x_{4} \partial_{2} + 2\,
         x_{2}  x_{5} x_{6} x_{7} \partial_{1}$\\
         &$ + 2\,
         x_{3}  x_{4} \partial_{3} + 2\,
          {{x_{1}}^{(2)}} \partial_{5}
    +  x_{1}  x_{2}\partial_{6} + 2\,
        x_{1} x_{3} \partial_{7} +  x_{1} x_{6} x_{7}\partial_{6} + 2\,
         {{x_{2}}^{(2)}}\partial_{7}
    + x_{2} x_{5}  x_{6}\partial_{5} + 2\,
          x_{2} x_{5}x_{7}\partial_{6}$\\
          &$ +  x_{3}  x_{5} x_{6}\partial_{6} + 2\,
      x_{3}  x_{5} x_{7}\partial_{7} +  x_{4} x_{5}\partial_{5}
    +  x_{4} x_{6} \partial_{6} +  x_{4} x_{7}\partial_7 $\\ \hline
\end{tabular}
$$
Let us study now the case where
$\fg_0'=\mathfrak{der}_0(\fg_{-})$. Our calculation shows that
$\fg_0'$ is generated by the vectors $Y_1, Z_1, H_1, H_2$ above
together with $ V=2 x_3 \partial_1+x_7 \partial_5$. The Lie
algebra $\fg_0'$ is solvable of $\sdim=5|0$. The CTS prolong
$(\fg_{-},\fg_0')_*$ gives a Lie superalgebra that is not simple
because $\fg_1'$ does not generate the positive part. Its simple
part is a new Lie superalgebra that we denote by $\mathfrak{BRJ}$,
described as follows (here also $N=(1,1,1)$:
$$
\renewcommand{\arraystretch}{1.4}
\begin{tabular}{|l|c|c|c|c|c|c|c|} \hline
&$\fg_1'$ & $\mathrm{ad}_{\fg_1'}(\fg_1')$&
$\mathrm{ad}_{\fg_1'}^2(\fg_1')$& $\mathrm{ad}_{\fg_1'}^3(\fg_1')$
& $\mathrm{ad}_{\fg_1'}^4(\fg_1')$
& $\mathrm{ad}_{\fg_1'}^5(\fg_1')$\\
\hline $\sdim$ &$0|6$ & $6|0$& $0|5$& $3|0$& $0|3$& $1|0$\\
\hline
\end{tabular}
$$

\ssec{A description of $\mathfrak{Brj}(3;\underline{N}|3)$} We have
the following realization of the non-positive part inside
$\fvect(3;\underline{N}|3)$:
\begin{equation}\label{tildeBrj(3|3)}
\small
\renewcommand{\arraystretch}{1.4}
\begin{tabular}{|l|l|} \hline
$\fg_{i}$&the generators (even $\mid$ odd) \\ \hline \hline
$\fg_{-2}$&$Y_7=\partial_1 \mid Y_5= \partial_4 ,\; Y_8=\partial_5$ \\ \hline
$\fg_{-1}$&$Y_{1}=\partial_2,\; Y_{6}=2\, x_2\partial_1+\partial_3\;\mid
Y_3=x_2\partial_4+x_3\partial_5+ 2\,x_6\partial_1+\partial_6$\\
\hline $\fg_{0}$&$H_2=[X_2,Y_2],\; H_1=x_1\partial_1+x_3\partial_3+ 2\,x_4\partial_4+
2\,x_6\partial_6, \; X_4=[X_2,X_2],\; Y_4=[Y_2,Y_2]\mid $\\
&$Y_2=x_2^{(2)}\partial_4+x_2\,
   x_3\partial_5+ 2\, x_2\,
   x_6\partial_1+ x_1\partial_5+ x_2\partial_6+ x_4\partial_1
   +x_6\partial_3$\\
    &$X_2=x_3^{(2)}\partial_5+ 2\, x_1\partial_4+x_3\partial_6+ x_5\partial_1+2\, x_6\partial_2$\\
\hline
\end{tabular}
\end{equation}
The Lie superalgebra $\fg_0$ is isomorphic to
$\mathfrak{osp}(1|2)\oplus \mathbb{K}$. The CTS prolong
$(\fg_{-},\fg_0)_*$ is NOT simple since it gives back
$\mathfrak{brj}(2;3)$ + an outer derivation. The $\sdim$ of the
positive parts are described as follows:
$$
\renewcommand{\arraystretch}{1.4}
\begin{tabular}{|l|c|c|c|} \hline
&$\fg_1$ & $\fg_2$& $\fg_3$\\
\hline $\sdim$ &$2|1$ & $1|2$& $0|1$\\
\hline
\end{tabular}
$$
and the lowest weight vectors are
$$ \tiny
\renewcommand{\arraystretch}{1.4}
\begin{tabular}{|l|l|}
\hline $V_1'$& $ 2\,x_1\,
   x_2\partial_1+x_2\,
   x_4\partial_4+x_3\,
   x_4\partial_5+2\, x_4\,
   x_6\partial_1+x_1\partial_3+2\,x_2\,
   x_3\partial_3+x_2\,
   x_6\partial_6+x_4\partial_6$\\\hline
$V_2'$ & $ x_1^{(2)}\partial_5+x_1\,
   x_2^{(2)}\partial_4+x_1\,
   x_2\, x_3\partial_5+2\,x_1\, x_2\,
   x_6\partial_1+x_1\,
   x_4\partial_1+2\,x_2^{(2)}\,
   x_3^{(2)}\partial_5+2\,x_2^{(2)}\,
   x_5\partial_1+2\, x_4\,
   x_5\partial_5+x_1\,
   x_2\partial_6$\\
   & $+x_1\,x_6\partial_3+2\,x_2^{(2)}\,
   x_3\partial_6+x_2^{(2)}\,
   x_6\partial_2+x_2\, x_3\,
   x_6\partial_3+x_2\,
   x_4\partial_2+x_2\,
   x_5\partial_3+2\, x_4\,
   x_6\partial_6$\\\hline
    $V_3'$ & $x_1^{(2)}\,   x_2\partial_4+x_1^{(2)}\,
   x_3\partial_5+2\, x_1^{(2)}\,
   x_6\partial_1+2\, x_1\,   x_2\, x_3^{(2)}\partial_5+2\,
   x_1\, x_2\, x_5\partial_1+x_2\, x_3^{(2)}\, x_4\partial_1+2\,
   x_2\, x_4\,
   x_5\partial_4+2\, x_3\,
   x_4\, x_5\partial_5$\\
   &$+x_4\,
   x_5\,   x_6\partial_1+x_1^{(2)}\partial_6+2\, x_1\, x_2\,
   x_3\partial_6+x_1\, x_2\,
   x_6\partial_2+x_1\, x_3\,
   x_6\partial_3+x_1\,x_4\partial_2+x_1\,
   x_5\partial_3+2\, x_2\,
   x_3^{(2)}\, x_6\partial_3+2\,x_2\, x_3\,
   x_4\partial_2$\\
   &$+2\, x_2\,
   x_3\, x_5\partial_3+x_2\,
   x_5\, x_6\partial_6+2\,x_3\, x_4\,
   x_6\partial_6+2\, x_4\,
   x_5\partial_6$\\ \hline
\end{tabular}
$$
\normalsize

\ssec{A description of $\mathfrak{Brj}(3;\underline{N}|4)$} We have
the following realization of the non-positive part inside
$\fvect(3;\underline{N}|4)$:
\begin{equation}\label{Brj(3|4)'}
\small
\renewcommand{\arraystretch}{1.4}
\begin{tabular}{|l|l|} \hline
$\fg_{i}$&the generators (even $\mid$ odd) \\ \hline \hline
$\fg_{-3}$&$\mid Y_{8}=\partial_4 $\\ \hline
$\fg_{-2}$&$Y_{4}=\partial_1,\; Y_{6}=\partial_2,\; Y_7=\partial_3\;\mid$\\
\hline $\fg_{-1}$&$\mid Y_2= x_{3}\partial_{4} +
2 x_5 \partial_{1}+\partial_5,\; Y_3= \partial_6+ x_5\,x_6\partial_4+
 x_2 \partial_4+x_5\partial_2+ 2\,x_6\partial_3$\\
&$Y_5=x_1 \partial_4 + x_5 \partial_3+ \partial_7$\\\hline
$\fg_{0}\simeq \mathfrak{hei}(2|0)\subplus \Kee H_2$&
$H_1=[Z_1,Y_1],\;H_2=2\, x_1\partial_1+x_2\partial_2+x_4\partial_4
+x_5\partial_5+ 2\,
   x_7\partial_7$\\
   & $ Y_1=x_1\partial_2+x_2\partial_3+ x_5\,
   x_6\partial_3+ 2\,
   x_5\partial_6+ 2\,
   x_6\partial_7$\\
   &$ Z_1= 2\, x_5\, x_6\,x_7\partial_4+2\,
   x_2\partial_1+x_3\partial_2+x_5\,
   x_6\partial_1+ 2\, x_6\,
   x_7\partial_3+ 2\,
   x_6\partial_5+x_7\partial_6\;\mid $
   \\
\hline
\end{tabular}
\end{equation}
The Lie superalgebra $\fg_0$ is solvable with the property that
$[\fg_0,\fg_0]=\mathfrak{hei}(2|0)$. The CTS prolong
$(\fg_{-},\fg_0)_*$ is NOT simple since $\fg_1$ does not generate
the positive part. Our calculation shows that the prolong does not
depend on $N$, i.e., $N=(1,1,1,1)$. The simple part of this prolong
is ${}^{3)}\mathfrak{brj}(2; 3)$. The $\sdim$ of the positive parts
are described as follows:
$$
\renewcommand{\arraystretch}{1.4}
\begin{tabular}{|l|c|c|c|} \hline
&$\fg_1$ & $\fg_2$& $\fg_3$\\
\hline $\sdim$ &$0|3$ & $3|0$& $0|2$\\
\hline
\end{tabular}
$$
and the lowest weight vectors are
$$ \tiny
\renewcommand{\arraystretch}{1.4}
\begin{tabular}{|l|l|}
\hline $V_1'$& $2 x_1\,
   x_3\partial_4+x_2^{(2)}\partial_4+x_2\, x_5\,
   x_6\partial_4+x_1\,
   x_5\partial_1+x_2\,
   x_5\partial_2+2\,x_2\,x_6\partial_3+2\,x_3\,
   x_5\partial_3+x_4\partial_3+2\,
   x_1\partial_5+x_2\partial_6+2\,
   x_3\partial_7+2\, x_5\,
   x_7\partial_7$\\\hline
$V_2'$ & $2\, x_1\,
   x_4\partial_4+x_1^{(2)}\partial_1+2\, x_1\,
   x_2\partial_2+2\,
   x_2^{(2)}\partial_3+2\,x_2\, x_5\,
   x_6\partial_3+2\,x_4\,
   x_5\partial_3+2\,x_1\,
   x_5\partial_5+x_1\,
   x_7\partial_7+x_2\,
   x_5\partial_6+x_2\,
   x_6\partial_7+2\,
   x_4\partial_7 $\\\hline
    $V_3'$ & $ x_1\,
   x_3^{(2)}\partial_4+2\,x_2^{(2)}\,
   x_3\partial_4+2\,x_2\, x_3\, x_5\,
   x_6\partial_4+2\,x_1\, x_3\,
   x_5\partial_1+2\,x_1\,
   x_4\partial_1+x_2^{(
   2)}\, x_5\partial_1+2\,x_2\, x_3\,
   x_5\partial_2+x_2\,
   x_3\, x_6\partial_3+2\,x_2\,
   x_4\partial_2$\\
   &$+2\,x_2\, x_5\, x_6\,
   x_7\partial_3+x_3^{(
   2)}\, x_5\partial_3+2\,x_3\,
   x_4\partial_3+x_1\,
   x_3\partial_5+x_1\,
   x_6\, x_7\partial_6+2\,
   x_2^{(2)}\partial_5+2\,x_2\,
   x_3\partial_6+2\,x_2\, x_5\,
   x_6\partial_5+x_2\,
   x_5\, x_7\partial_6$\\
   &$+2\,x_2\, x_6\,
   x_7\partial_7+x_3^{(2)}\partial_7+x_3\, x_5\,
   x_7\partial_7+x_4\,
   x_5\partial_5+x_4\,
   x_6\partial_6+x_4\,
   x_7\partial_7 $\\ \hline
\end{tabular}
$$
Let us study now the case where $\fg_0'=\mathfrak{der}_0(\fg_{-})$.
The Lie algebra $\fg_0'$ is solvable of $\sdim=5|0$. The CTS prolong
$(\fg_{-},\fg_0')_*$ gives a Lie superalgebra that is not simple
because $\fg_1'$ does not generate the positive part. Its simple
part is a new Lie superalgebra that we had denoted by
$\mathfrak{BRJ}$, described as follows (here also $N=(1,1,1)$:
$$
\renewcommand{\arraystretch}{1.4}
\begin{tabular}{|l|c|c|c|c|c|c|c|} \hline
&$\fg_1'$ & $\mathrm{ad}_{\fg_1'}(\fg_1')$&
$\mathrm{ad}_{\fg_1'}^2(\fg_1')$& $\mathrm{ad}_{\fg_1'}^3(\fg_1')$
& $\mathrm{ad}_{\fg_1'}^4(\fg_1')$
& $\mathrm{ad}_{\fg_1'}^5(\fg_1')$\\
\hline $\sdim$ &$0|6$ & $6|0$& $0|5$& $3|0$& $0|3$& $1|0$\\
\hline
\end{tabular}
$$

\ssec{Constructing Melikyan superalgebras} Denote by
$\cal{F}_{1/2}:=\cO(1;1)\sqrt{dx}$ the space of semi-densities
(weighted densities of weight $\frac12$). For $p=3$, the CTS prolong
of the triple $(\mathbb{K}, \Pi(\cal{F}_{1/2}), \fcvect(1;1))_*$
gives the whole $\fk(1;\underline{N}|3)$. For $p=5$, let us realize
the non-positive part in $\fk(1;\underline{N}|5)$:
\begin{equation}\label{???}
\small
\renewcommand{\arraystretch}{1.4}
\begin{tabular}{|l|l|} \hline
$\fg_{i}$&the generators  \\ \hline \hline
$\fg_{-2}$&$1$\\
\hline $\fg_{-1}$&$\Pi(\cal{F}_{1/2})$\\
   \hline $\fg_{0}$& $\partial_1\longleftrightarrow 4\,  \xi_{1} \eta_2+  \xi_2 \theta,\;
   x_1\partial_1\longleftrightarrow 2\,  \xi_{1}  \eta_1 +\xi_2  \eta_2,\;
   x_1^{(2)}\partial_1\longleftrightarrow 2\,  \xi_{2}  \eta_1 + 3\,  \theta\eta_2
   ,\;  x_1^{(3)}\partial_1\longleftrightarrow 2\, \theta \eta_1$\\
   &$ x_1^{(4)}\partial_1\longleftrightarrow 2\,  \eta_2\eta_1,\qquad t  $
   \\
\hline
\end{tabular}
\end{equation}
The CTS prolong gives that $\fg_i$=0 for all $i>0$.

Consider now the case of  $(\mathbb{K}, \Pi(\cal{F}_{1/2}),
\fcvect(2;1))_*$, where $p=3$. The non-positive part is realized in
$\fk(1;\underline{N}|9)$ as follows:
\begin{equation}\label{fk(1;N|9)}
\small
\renewcommand{\arraystretch}{1.4}
\begin{tabular}{|l|l|} \hline
$\fg_{i}$&the generators  \\ \hline \hline
$\fg_{-2}$&$1$\\
\hline $\fg_{-1}$&$\Pi(\cal{F}_{1/2})$\\
   \hline $\fg_{0}$& $\partial_1 \longleftrightarrow 2\, \xi_{1} \eta_3
   + x_{2} \theta + 2\,  \xi_3 \eta_4,\; x_1\partial_1 \longleftrightarrow  \xi_{1}\eta_1
   + \xi_{2} \eta_2 + 2\, \xi_4  \eta_4,\;x_1^2\partial_1 \longleftrightarrow  \xi_3 \eta_1
   + \xi_4  \eta_3
   + \theta \eta_2,$\\
   &$\;\partial_2 \longleftrightarrow 2\,  \xi_{1}  \eta_2
   + \xi_{2} \xi_4
   + \xi_3  \theta ,\;x_2 \partial_2\longleftrightarrow  \xi_{1} \eta_1
   +  \xi_3  \eta_3
   + \xi_4 \eta_4,\;x_2^2 \partial_2\longleftrightarrow \xi_{2} \eta_1
   +  \theta  \eta_3
   + 2\, \eta_4\eta_2,$\\
   &$x_1 x_2 \partial_1\longleftrightarrow \xi_{2}\eta_1
   + \eta_4 \eta_2 ,\;x_1 x_2 \partial_2\longleftrightarrow  \xi_3  \eta_1
   + 2\, \xi_4 \eta_3,\;x_1^2 x_2 \partial_1\longleftrightarrow  \theta  \eta_1
   + 2\, \eta_3 \eta_2,$\\
   &$x_1^2 x_2 \partial_2\longleftrightarrow \xi_4  \eta_1 ,\;
   x_1x_2^2 \partial_1\longleftrightarrow 2\, \eta_4  \eta_1 ,\;x_1x_2^2 \partial_2
   \longleftrightarrow \theta  \eta_1
   + \eta_3 \eta_2 ,\; x_1^2x_2^2 \partial_1\longleftrightarrow \eta_3 \eta_1
   ,$\\
&$    x_1^2x_2^2 \partial_2\longleftrightarrow \eta_2\eta_1 ,\qquad
t $
   \\
\hline
\end{tabular}
\end{equation}
The CTS prolong $(\fg_{-},\fg_0)_*$ gives a Lie superalgebra that is
not simple with the property that $\sdim(\fg_1)=0|4$ and $\fg_i=0$
for all $i>1$. The generating functions  of $\fg_1$ are
$$
\xi_{2} \eta_2\eta_{1}
    + 2\,  \xi_{3} \eta_3 \eta_{1}
    + \xi_4  \eta_4 \eta_{1}
    + \theta \eta_3 \eta_2 ,\quad
  2\, \xi_4 \eta_3 \eta_{1}
    + \theta  \eta_2  \eta_{1}  ,\quad
  \theta \eta_3 \eta_{1}
    + \eta_4  \eta_2 \eta_{1} ,\quad
\eta_3  \eta_2  \eta_{1}.
$$

\ssec{Defining relations of the positive parts of
$\mathfrak{brj}(2;3)$ and $\mathfrak{brj}(2;5)$} For the
presentations of the Lie superalgebras with Cartan matrix, see
\cite{GL1, BGL1}. The only nontrivial part of these relations are
analogs of the Serre relations (both the straightforward ones and
the ones different
in shape). Here they are:\\

\noindent \underline{$\mathfrak{brj}(2;3)$};\quad$\sdim\;
\mathfrak{brj}(2;3)=10|8$.

$ \arraycolsep=0pt
\renewcommand{\arraystretch}{1.2}
\begin{array}{ll}
1)\ &{{\left[\left[x_{1},\,x_{2}\right],\,
      \left[x_{2},\,\left[x_{1},\,x_{2}\right]\right]\right]}=0},\\
&      {{\left[\left[x_{2},\,x_{2}\right],\,
      \left[\left[x_{1},\,x_{2}\right],\,
       \left[x_{2},\,x_{2}\right]\right]\right]}=0}.
      \end{array}
$

$ \arraycolsep=0pt
\renewcommand{\arraystretch}{1.2}
\begin{array}{ll}
2)\ & \ad_{x_{2}}^3(x_{1})=0,\\
&  {{\left[\left[x_{1},\,x_{2}\right],\,
      \left[\left[x_{1},\,x_{2}\right],\,
       \left[x_{1},\,x_{2}\right]\right]\right]}=0},\\
&  {{\left[\left[x_{2},\,\left[x_{1},\,x_{2}\right]\right],\,
      \left[\left[x_{1},\,x_{2}\right],\,
       \left[x_{2},\,\left[x_{1},\,x_{2}\right]\right]\right]
       \right]}=0}.
      \end{array}
$

$ \arraycolsep=0pt
\renewcommand{\arraystretch}{1.2}
\begin{array}{ll}
3)\ & \ad_{x_{1}}^3(x_{2})=0,\\
&    [x_2,[x_1,[x_1,x_2]]]-
   [[x_1,x_2],[x_1,x_2]]=0,\\
 &  [[x_1,x_2],[x_2,x_2]]=0.
      \end{array}
$

\noindent \underline{$\mathfrak{brj}(2;5)$}; \quad$\sdim\;
\mathfrak{brj}(2;5)=10|12$.

$ \arraycolsep=0pt
\renewcommand{\arraystretch}{1.2}
\begin{array}{ll}
1)\ &  {{\left[\left[x_{2},\,\left[x_{1},\,x_{2}\right]\right],\,
      \left[x_{2},\,\left[x_{1},\,x_{2}\right]\right]\right]}
={2\, \left[\left[x_{1},\,x_{2}\right],\,
        \left[\left[x_{1},\,x_{2}\right],\,
         \left[x_{2},\,x_{2}\right]\right]\right]}},\\
&  {{\left[\left[x_{2},\,x_{2}\right],\,
      \left[\left[x_{1},\,x_{2}\right],\,
       \left[x_{2},\,x_{2}\right]\right]\right]}=0},\\
&  {{\left[\left[x_{2},\,\left[x_{1},\,x_{2}\right]\right],\,
      \left[\left[x_{1},\,x_{2}\right],\,
       \left[x_{2},\,\left[x_{1},\,x_{2}\right]\right]\right]
       \right]}= 0}.
      \end{array}
$

$ \arraycolsep=0pt
\renewcommand{\arraystretch}{1.2}
\begin{array}{ll}
2\ )&  \ad_{x_{2}}^4(x_{1})=0,\\
&  {{\left[\left[x_{2},\,\left[x_{1},\,x_{2}\right]\right],\,
      \left[x_{2},\,\left[x_{2},\,\left[x_{1},\,x_{2}\right]
        \right]\right]\right]}=0},\\
&  {{\left[\left[\left[x_{1},\,x_{2}\right],\,
       \left[x_{1},\,x_{2}\right]\right],\,
      \left[\left[x_{1},\,x_{2}\right],\,
       \left[x_{2},\,\left[x_{1},\,x_{2}\right]\right]\right]
       \right]}=0}.
      \end{array}
$

\end{document}